\documentclass[10pt]{amsart} 

\usepackage{amsmath}
\usepackage{amssymb}

\usepackage{graphicx}
\usepackage{xcolor}
\usepackage{subcaption}
\usepackage{multirow}
\usepackage{enumitem}
\usepackage{mathrsfs}
\usepackage[rightcaption]{sidecap}
\usepackage{rotating}
\usepackage{tikz}

\newtheorem{theorem}{Theorem}[section]

\theoremstyle{definition}
\newtheorem{example}[theorem]{Example}
\theoremstyle{remark}

\newcommand{\R}{\mathbb R}
\newcommand{\bfu}{\mathbf u}
\newcommand{\bfx}{\mathbf x}
\newcommand{\bff}{\mathbf f}
\newcommand{\bfg}{\mathbf g}
\newcommand{\bfs}{\mathbf s}
\newcommand{\dx}{\mathrm d}

\newcommand{\JS}{\mathrm{JS}}
\newcommand{\Z}{\mathrm Z}

\newcommand{\CADNN}{\mathrm{CADNN}}
\newcommand{\loss}{\mathcal L}
\newcommand{\CAD}{\mathrm{CAD}}

\newcommand{\SYM}{\mathrm{SYM}}
\newcommand{\LN}{\mathrm{LN}}

\newcommand{\e}{\mathrm e}
\newcommand{\mdf}{\mathfrak m}

\newcommand{\cfl}{\mathrm{CFL}}
\newcommand{\weno}{\mathcal{WENO}}
\newcommand{\N}{\mathcal N}

\newcommand{\vecS}{\mathcal{S}}

\title[Conservative approximation-based FNN for WENO schemes]{Conservative approximation-based feedforward neural network for WENO schemes}
\author{Kwanghyuk Park}
\address{Graduate School of Artificial Intelligence $\&$ POSTECH MINDS (Mathematical Institute for Data Science), Pohang University of Science and Technology, Pohang 37673, Korea}
\email{pkh0219@postech.ac.kr}

\author{Jiaxi Gu}
\address{Department of Mathematics $\&$ POSTECH MINDS (Mathematical Institute for Data Science), Pohang University of Science and Technology, Pohang 37673, Korea}
\email{jiaxigu@postech.ac.kr}

\author{Jae-Hun Jung}
\address{Department of Mathematics $\&$ POSTECH MINDS (Mathematical Institute for Data Science), Pohang University of Science and Technology, Pohang 37673, Korea}
\email{jung153@postech.ac.kr}

\makeatletter
\@namedef{subjclassname@2020}{\textup{}2020 Mathematics Subject Classification}
\makeatother

\begin{document}

\maketitle

\begin{abstract}
In this work, we present the feedforward neural network based on the conservative approximation to the derivative from point values, for the weighted essentially non-oscillatory (WENO) schemes in solving hyperbolic conservation laws. 
The feedforward neural network, whose inputs are point values from the three-point stencil and outputs are two nonlinear weights, takes the place of the classical WENO weighting procedure. 
For the training phase, we employ the supervised learning and create a new labeled dataset for one-dimensional conservative approximation, where we construct a numerical flux function from the given point values such that the flux difference approximates the derivative to high-order accuracy.
The symmetric-balancing term is introduced for the loss function so that it propels the neural network to match the conservative approximation to the derivative and satisfy the symmetric property that WENO3-JS and WENO3-Z have in common.
The consequent WENO schemes, WENO3-CADNNs, demonstrate robust generalization across various benchmark scenarios and resolutions, where they outperform WENO3-Z and achieve accuracy comparable to WENO5-JS.
\end{abstract}

\section{Introduction} \label{sec:intro}
We consider the system of hyperbolic conservation laws,
\begin{equation} \label{eq:hyperbolic}
 \bfu_t + \sum_{i=1}^d \bff_i(\bfu)_{x_i} = 0,
\end{equation}
where $\bfu(\bfx,t) = (u_1(\bfx,t), \cdots, u_m(\bfx,t))$ is the vector of $m$ conserved quantities with $\bfx = (x_1, \cdots, x_d)$ and $d$ the spatial dimension, and each $\bff_i(\bfu)$ is the flux function in the spatial direction $x_i$.
The system is supplemented with the initial condition $\bfu(\bfx,0)=\bfu_0(\bfx)$ and appropriate boundary conditions.
As the solution may develop discontinuities, robust numerical schemes are required to accurately capture the features while suppressing spurious oscillations arising from the Gibbs phenomenon.
The essentially non-oscillatory (ENO) and weighted ENO (WENO) schemes have been proven particularly effective in resolving discontinuities and minimizing nonphysical oscillations near large gradients.

Harten \cite{Harten83} initiated a class of explicit second-order total-variation-diminishing (TVD) finite difference schemes.
To achieve higher-order accuracy, Harten et al. \cite{Harten} introduced ENO scheme within the finite volume framework.
Later in \cite{ShuOsherI}, Shu and Osher developed the finite difference ENO scheme using the fluxes to simplify the reconstruction procedure, particularly for multi-dimensional problems.
Jiang and Shu \cite{Jiang} proposed the fifth-order finite difference WENO scheme (WENO5-JS) with the formulation of smoothness indicators.
In order to reduce the numerical dissipation around discontinuities, Borges et al. \cite{Borges} first constructed the global smoothness indicator, and presented the Z-type nonlinear weights $\omega^{\Z}_k$ and WENO5-Z scheme.
These fifth-order WENO schemes can be easily adapted to their third-order analogs, WENO3-JS and WENO3-Z \cite{Shu,Don}.

In recent developments, machine learning has increasingly been applied to the WENO schemes, unlocking the potential for the improvements in both accuracy and efficiency.
In \cite{Sun,Wen,Xue}, the neural network-based discontinuity detectors were applied to the hybrid WENO schemes.
Kossaczk\'{a} et al. \cite{Kossaczka} admitted a small convolutional neural network, adjusting smoothness indicators to the fifth-order WENO scheme.
In \cite{Wang}, the reinforcement learning policy network was employed to optimize the nonlinear weights of the numerical fluxes corresponding to the substencils.
Bezgin et al. \cite{Bezgin21} designed a convolutional neural network to produce the nonlinear weights and the dispersion coefficient at the cell face of the finite volume scheme for the linear diffusive-dispersive regularizations of the scalar cubic conservation law.
Subsequently in \cite{Bezgin22}, the Delta layer, which pre-processes the input data, was introduced to the multilayer perceptron, predicting the nonlinear weights for the third-order finite volume WENO scheme.
In addition, for the simulation, the output weights would pass through an ENO layer in order to retain the ENO property.
In \cite{Nogueira24}, Nogueira et al. presented a neural network including the ENO layer, to calculate the nonlinear weights for the fifth-order finite difference WENO scheme.
For computational efficiency, we \cite{Park} proposed the shallow neural network without the ENO layer, serving as a WENO weighting function for the third-order finite difference WENO scheme. 

In this work, we use the feedforward neural network as a weighting function in WENO schemes.
The architecture of the neural network is composed of the input layer, the pre-processing layer, two hidden layers, and the output layer, where we modify the Delta layer in \cite{Bezgin22,Park} as the pre-processing layer for stability and skip the ENO layer in \cite{Bezgin22,Nogueira24} for computational efficiency.
We then create a new labeled dataset for the one-dimensional conservative approximation to the derivative from point values with the uniform grid, where a numerical flux function is constructed so that the flux difference approximates the derivative to high-order accuracy \cite{Shu}.
During the training phase, we employ the supervised learning to optimize two nonlinear weights for different three-point stencils.
In order to minimize the difference between the predicted derivative and the label, we introduce the error of the conservative approximation to derivative for the loss function. 
The symmetric-balancing term is designed for the loss function, which propels the neural network to follow the symmetric property that we observe in the classical WENO3-JS and WENO3-Z schemes.
We also include the linear term \cite{Bezgin22,Park} in the loss function to reduce numerical dissipation around discontinuities and limit the nonlinear weights to linear weights in smooth regions \cite{Gu,ChenGuJung}.
Numerical simulations in both one and two space dimensions demonstrate the improved performance of the proposed WENO3-CADNN schemes over WENO3-Z for all benchmark problems and the comparable resolution to WENO5-JS in some tests.

The structure of this paper is as follows.
We start in Section \ref{sec:weno} with a brief review of the third-order finite difference WENO schemes.
Section \ref{sec:fnn} details the architecture of the feedforward neural network, the labeled dataset for the supervised learning, and the training process.
The numerical results for one- and two-dimensional test problems are presented in Section \ref{sec:nr}.
We end in Section \ref{sec:conclusions} with the conclusion.

\section{Third-order WENO schemes} \label{sec:weno}
In this section, we briefly describe the third-order WENO schemes in the finite difference framework. 
Consider the one-dimensional scalar hyperbolic conservation law,
\begin{equation} \label{eq:1d_scalar_hyperbolic}
 \frac{\partial}{\partial t} u(x,t) + \frac{\partial}{\partial x} f(u(x,t)) = 0.
\end{equation}
We discretize the spatial domain into $N$ cells of the uniform size $\Delta x$ where $x_i$ is the cell center of the $i$th cell $I_i = (x_{i-1/2},x_{i+1/2})$ and $x_{i\pm1/2} = x_i \pm \Delta x/2$ are the cell faces. 
If we fix the time $t$, the semi-discretization of \eqref{eq:1d_scalar_hyperbolic} gives 
\begin{equation} \label{eq:1D_hyperbolic_discrete}
 \frac{du(x_i, t)}{dt} = - \left. \frac{\partial f \left( u(x,t) \right)}{\partial x} \right |_{x=x_i}.
\end{equation} 
Then we define the auxiliary function $h(x)$ implicitly by
$$
   f \left( u(x) \right) = \frac{1}{\Delta x} \int^{x+\Delta x/2}_{x-\Delta x/2} h(\xi) \dx \xi, 
$$
where the time variable $t$ is dropped to simplify the notation. 
Differentiating both sides with respect to $x$, we obtain
\begin{equation} \label{eq:partial_f}
 \frac{\partial f}{\partial x} = \frac{h(x+\Delta x/2) - h(x-\Delta x/2)}{\Delta x}.
\end{equation}
With the evaluation of \eqref{eq:partial_f} at $x=x_i$, \eqref{eq:1D_hyperbolic_discrete} becomes 
\begin{equation} \label{eq:1D_hyperbolic_discrete_h}
 \frac{du(x_i,t)}{dt} = - \frac{h_{i+1/2} - h_{i-1/2}}{\Delta x},
\end{equation} 
with $h_{i \pm 1/2} = h(x_{i \pm 1/2})$.
We aim to reconstruct the fluxes $h_{i \pm 1/2}$ at the cell faces such that the numerical approximation achieves high-order accuracy in smooth regions and inhibits oscillations around discontinuities.

\begin{figure}[htbp]
\vspace{-0.15in}
\centering
\begin{tikzpicture}
 \draw (-3,0) -- (3,0);
 \filldraw[red] (-3,0) circle (2pt);
 \filldraw[red] (0,0) circle (2pt);
 \filldraw[red] (3,0) circle (2pt);
 \node[rectangle, inner sep=2pt, minimum size=2.5pt, draw=blue] (12) at (1.5,0) {};
 \draw [->,thick] (1.5,0.7)--(12);
 \draw [blue,->,thick] (0.5,0.4)--(1.5,0.4);
 \node at (0.95,0.67) {\scriptsize $\textcolor{blue}{\hat{f}^+_{i+1/2}}$};
 \node at (-3+0.07,-1/3) {\scriptsize $f_{i-1}^{+}$};
 \node at (0,-1/3) {\scriptsize $f_{i}^{+}$};
 \node at (3+0.07,-1/3) {\scriptsize $f_{i+1}^{+}$};
 \node at (-3.3,-0.7) {\textcolor{black}{\scriptsize $\vecS$}};
 \draw (-3,-3/4) -- (3,-3/4);
 \filldraw[blue] (-3,-3/4) circle (2pt);
 \filldraw[blue] (0,-3/4) circle (2pt);
 \filldraw[blue] (3,-3/4) circle (2pt);
 \node at (3.3,-3/4) {\textcolor{black}{\scriptsize $\tau_3$}};
 \node at (-3.3,-13/12) {\scriptsize$S_0$};
 \draw (-3,-13/12) -- (0,-13/12);
 \filldraw[black] (-3,-13/12) circle (2pt);
 \filldraw[black] (0,-13/12) circle (2pt);
 \node at (0.3,-13/12) {\scriptsize$\beta_0$};
 \node at (-0.3,-17/12) {\scriptsize $S_1$};
 \draw (0,-17/12) -- (3,-17/12);
 \filldraw[black] (0,-17/12) circle (2pt);
 \filldraw[black] (3,-17/12) circle (2pt);
 \node at (3.3,-17/12) {\scriptsize $\beta_1$};
\end{tikzpicture}
\caption{The construction of the numerical flux $\hat{f}^+_{i+1/2}$ depends on the stencil $\left( f^+_{i-1}, f^+_i, f^+_{i+1} \right)$.}
\label{fig:stencil_l}
\end{figure}
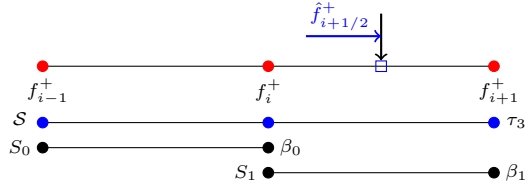

\begin{figure}[htbp]
\vspace{-0.15in}
\centering
\begin{tikzpicture}
 \draw (-3,0) -- (3,0);
 \filldraw[red] (-3,0) circle (2pt);
 \filldraw[red] (0,0) circle (2pt);
 \filldraw[red] (3,0) circle (2pt);
 \node[rectangle, inner sep=2pt, minimum size=2pt, draw=blue] (12) at (-1.5,0) {}; 
 \draw [->,thick] (-1.5,0.7)--(12);
 \draw [blue,->,thick] (-0.5,0.4)--(-1.5,0.4);
 \node at (-0.95,0.67) {\scriptsize $\textcolor{blue}{\hat{f}^-_{i+1/2}}$};
 \node at (-3+0.03,-1/3) {\scriptsize $f_i^{-}$};
 \node at (0+0.1,-1/3) {\scriptsize $f_{i+1}^{-}$};
 \node at (3+0.1,-1/3) {\scriptsize $f_{i+2}^{-}$};
 \node at (-3.3,-0.7) {\textcolor{black}{\scriptsize $\vecS$}};
 \draw (-3,-3/4) -- (3,-3/4);
 \filldraw[blue] (-3,-3/4) circle (2pt);
 \filldraw[blue] (0,-3/4) circle (2pt);
 \filldraw[blue] (3,-3/4) circle (2pt);
 \node at (3.3,-3/4) {\scriptsize\textcolor{black}{$\tau_3$}};
 \node at (-0.3,-13/12) {\scriptsize$S_0$};
 \draw (0,-13/12) -- (3,-13/12);
 \filldraw[black] (0,-13/12) circle (2pt);
 \filldraw[black] (3,-13/12) circle (2pt);
 \node at (3.3,-13/12) {\scriptsize $\beta_0$};
 \node at (-3.3,-17/12) {\scriptsize$S_1$};
 \draw (-3,-17/12) -- (0,-17/12);
 \filldraw[black] (-3,-17/12) circle (2pt);
 \filldraw[black] (0,-17/12) circle (2pt);
 \node at (0.4,-17/12) {\scriptsize$\beta_1$};
\end{tikzpicture}
\caption{The construction of the numerical flux $\hat{f}^-_{i+1/2}$ depends on the stencil $\left( f^-_{i+2}, f^-_{i+1}, f^-_i \right)$.}
\label{fig:stencil_r}
\end{figure}
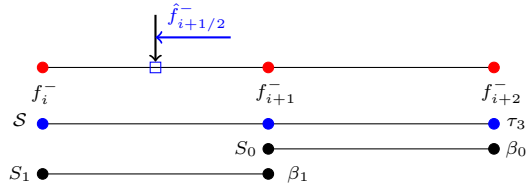

To approximate the flux $h_{i+1/2}$, we first use the Lax-Friedrichs flux splitting:
$$ f^{\pm}(u) = \frac{1}{2} \left( f(u) \pm \alpha u \right), $$
with $\alpha = \max_u |f'(u)|$ over the relevant range of $u$, and then construct the positive and negative numerical fluxes $\hat{f}^{\pm}_{i+1/2}$.
For the third-order WENO reconstruction, we take the stencil $\left( f^+_{i-1}, f^+_i, f^+_{i+1} \right)$ with $f^+_j = f^+(u_j), j = i-1,i,i+1$, for the positive numerical flux $\hat{f}^+_{i+1/2}$, as shown in Fig. \ref{fig:stencil_l}.
The numerical flux $\hat{f}^+_{i+1/2}$ takes a convex combination of two candidate numerical fluxes $\hat{f}^{k+}_{i+1/2}$ ($k=0,1$),
\begin{equation} \label{eq:positive_numerical_flux}
 \hat{f}^+_{i+1/2} = \omega_0 \hat{f}^{0+}_{i+1/2} + \omega_1 \hat{f}^{1+}_{i+1/2},
\end{equation}
where
$$ \hat{f}^{0+}_{i+1/2} = -\frac{1}{2} f^+_{i-1} + \frac{3}{2} f^+_i, \quad \hat{f}^{1+}_{i+1/2} = \frac{1}{2} f^+_i + \frac{1}{2} f^+_{i+1}. $$
In \cite{Shu}, the smoothness indicator $\beta_k$ of the substencil $S_k$ is defined as
$$
 \beta_0 = \left( f^+_{i-1} - f^+_i \right)^2, \quad \beta_1 = \left( f^+_i - f^+_{i+1} \right)^2.
$$
The nonlinear weights $\omega^{\JS}_k$ in \eqref{eq:positive_numerical_flux} are given by
$$
 \omega^{\JS}_k = \frac{\alpha_k}{\alpha_0 + \alpha_1}, \quad \alpha_k = \frac{d_k}{(\beta_k + \varepsilon)^2}, \quad k=0,1,
$$
where $d_0 = \frac{1}{3}, \, d_1 = \frac{2}{3}$ are linear weights and $\varepsilon$ is a small constant to prevent the denominator from being zero with $\varepsilon=10^{-6}$ in \cite{Jiang}.
With a different approach, Don and Borges \cite{Don} introduced the global smoothness indicator $\tau_3$,
$$ \tau_3 = \left| \beta_0 - \beta_1 \right|. $$
The Z-type nonlinear weights $\omega^\Z_k$ are defined as
$$
 \omega^\Z_k = \frac{\alpha_k}{\alpha_0 + \alpha_1}, \quad \alpha_k = d_k \left( 1 + \left( \frac{\tau_3}{\beta_k + \varepsilon} \right)^2 \right), \quad k=0,1,
$$
with $\varepsilon = 10^{-40}$ in \cite{Borges}.
The negative numerical flux $\hat{f}^-_{i+1/2}$, in Fig. \ref{fig:stencil_r}, is constructed in the mirror-symmetric way
$$ \hat{f}^-_{i+1/2} = \omega_0 \hat{f}^{0-}_{i+1/2} + \omega_1 \hat{f}^{1-}_{i+1/2}, $$
where
$$ \hat{f}^{0-}_{i+1/2} = -\frac{1}{2} f^-_{i+2} + \frac{3}{2} f^-_{i+1}, \quad \hat{f}^{1-}_{i+1/2} = \frac{1}{2} f^-_{i+1} + \frac{1}{2} f^-_i. $$
The nonlinear weights $\omega_k$ can be obtained using the same WENO weighting strategies as above.
The numerical flux $\hat{f}_{i+1/2}$, which is the sum of $\hat{f}^+_{i+1/2}$ and $\hat{f}^-_{i+1/2}$, is the approximation to $h_{i+1/2}$.

Note that each WENO weighting strategy can be regarded as a function $\weno^\N:\R^3 \to [0, \, 1]^2$ defined by 
$$ \weno^\N (\vecS)= \omega,~\N \in \{ \JS, \Z \}, $$ 
where $\vecS = (f_0, f_1, f_2)$ is the three-point stencil and $\omega = (\omega_0, \omega_1)$ is the set of two nonlinear weights.
We showed in \cite{Park} that those WENO weighting functions are scale-invariant if $\beta_0, \beta_1 \gg \varepsilon$ and translation-invariant.
While the classical WENO weighting functions are heuristic, we proposed the data-driven WENO weighting function using the shallow neural network \cite{Park} as an alternative.
The architecture of the shallow neural network includes the input layer, a pre-processing layer, one hidden layer and the output layer.
The pre-processing layer, called the Delta layer \cite{Bezgin22,Park}, transforms the three-point stencil $\vecS$ into four features $\Delta_j$ defined as
\begin{equation} \label{eq:delta_layer}
\begin{aligned}
 & \tilde{\Delta}_1 = \left| f_0-f_1 \right|, \, 
   \tilde{\Delta}_2 = \left| f_1-f_2 \right|, \,
   \tilde{\Delta}_3 = \left| f_0-f_2 \right|, \, 
   \tilde{\Delta}_4 = \left| f_0-2f_1+f_2 \right|, \\
 & \Delta_j = \tilde{\Delta}_j / \max \left( \tilde{\Delta}_1, \tilde{\Delta}_2, \varepsilon \right), \, \varepsilon = 10^{-12}, \, j=1,2,3,4.
\end{aligned}     
\end{equation}
In the training phase, we adopted the supervised learning, in which the training dataset of the three-point stencils is from the piecewise smooth function, and we computed the nonlinear weights $\omega^\JS_k$ by the WENO3-JS weighting function as labels.

\section{Feedforward neural network for WENO weighting function} \label{sec:fnn} 
We continue to construct the WENO weighting function in the data-driven approach using the feedforward neural network.  
We first describe the architecture in Section \ref{sec:fnn_architecture}. 
Afterward, we present the conservative approximation to the derivative from point values, where a new labeled dataset is created in Section \ref{sec:fnn_dataset}.
The training process is thoroughly discussed in Section \ref{sec:fnn_training}. 

\subsection{Architecture} \label{sec:fnn_architecture}
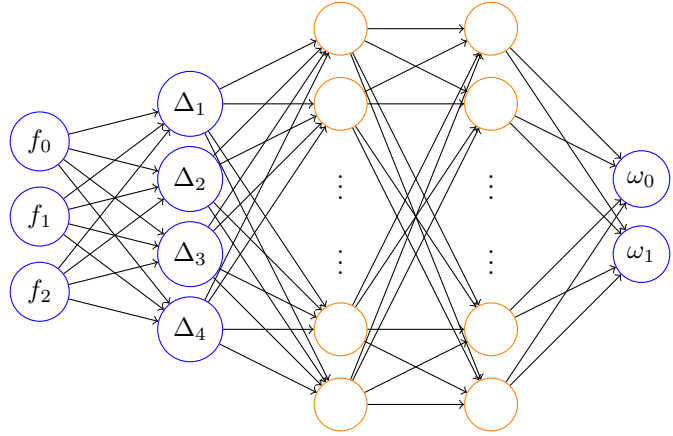
\begin{figure}[htbp]
\centering
\begin{tikzpicture}
\node[circle, minimum width=15pt, minimum height=15pt, draw=blue] (11) at(0,1)  {$f_0$};
\node[circle, minimum width=15pt, minimum height=15pt, draw=blue] (12) at(0,0)  {$f_1$};
\node[circle, minimum width=15pt, minimum height=15pt, draw=blue] (13) at(0,-1) {$f_2$};
\node[circle, minimum width=15pt, minimum height=15pt, draw=blue] (21) at(2,1.5)  {$\Delta_1$};
\node[circle, minimum width=15pt, minimum height=15pt, draw=blue] (22) at(2,0.5)  {$\Delta_2$};
\node[circle, minimum width=15pt, minimum height=15pt, draw=blue] (23) at(2,-0.5)  {$\Delta_3$};
\node[circle, minimum width=15pt, minimum height=15pt, draw=blue] (24) at(2,-1.5) {$\Delta_4$};
\node[circle, minimum width=15pt, minimum height=20pt, draw=orange] (31) at(4,2.5)  {};
\node[circle, minimum width=15pt, minimum height=20pt, draw=orange] (32) at(4,1.5)  {};
\node (33) at(4,0.5)  {$\vdots$};
\node (34) at(4,-0.5) {$\vdots$};
\node[circle, minimum width=15pt, minimum height=20pt, draw=orange] (35) at(4,-1.5) {};
\node[circle, minimum width=15pt, minimum height=20pt, draw=orange] (36) at(4,-2.5) {};
\node[circle, minimum width=15pt, minimum height=20pt, draw=orange] (41) at(6,2.5)  {};
\node[circle, minimum width=15pt, minimum height=20pt, draw=orange] (42) at(6,1.5)  {};
\node (43) at(6,0.5)  {$\vdots$};
\node (44) at(6,-0.5) {$\vdots$};
\node[circle, minimum width=15pt, minimum height=20pt, draw=orange] (45) at(6,-1.5)  {};
\node[circle, minimum width=15pt, minimum height=20pt, draw=orange] (46) at(6,-2.5)  {};
\node[circle, minimum width=15pt, minimum height=15pt, draw=blue] (51) at(8,0.5)  {$\omega_0$};
\node[circle, minimum width=15pt, minimum height=15pt, draw=blue] (52) at(8,-0.5)  {$\omega_1$};
\draw[->] (11)--(21);
\draw[->] (11)--(22);
\draw[->] (11)--(23);
\draw[->] (11)--(24);
\draw[->] (12)--(21);
\draw[->] (12)--(22);
\draw[->] (12)--(23);
\draw[->] (12)--(24);
\draw[->] (13)--(21);
\draw[->] (13)--(22);
\draw[->] (13)--(23);
\draw[->] (13)--(24);
\draw[->] (21)--(31);
\draw[->] (21)--(32);
\draw[->] (21)--(35);
\draw[->] (21)--(36);
\draw[->] (22)--(31);
\draw[->] (22)--(32);
\draw[->] (22)--(35);
\draw[->] (22)--(36);
\draw[->] (23)--(31);
\draw[->] (23)--(32);
\draw[->] (23)--(35);
\draw[->] (23)--(36);
\draw[->] (24)--(31);
\draw[->] (24)--(32);
\draw[->] (24)--(35);
\draw[->] (24)--(36);
\draw[->] (31)--(41);
\draw[->] (32)--(41);
\draw[->] (35)--(41);
\draw[->] (36)--(41);
\draw[->] (31)--(42);
\draw[->] (32)--(42);
\draw[->] (35)--(42);
\draw[->] (36)--(42);
\draw[->] (31)--(45);
\draw[->] (32)--(45);
\draw[->] (35)--(45);
\draw[->] (36)--(45);
\draw[->] (31)--(46);
\draw[->] (32)--(46);
\draw[->] (35)--(46);
\draw[->] (36)--(46);
\draw[->] (41)--(51);
\draw[->] (42)--(51);
\draw[->] (45)--(51);
\draw[->] (46)--(51);
\draw[->] (41)--(52);
\draw[->] (42)--(52);
\draw[->] (45)--(52);
\draw[->] (46)--(52);
\end{tikzpicture}
\caption{The feedforward neural network architecture.}
\label{fig:dnn}
\end{figure}

The feedforward neural network is structured with the input layer, a pre-processing layer, two hidden layers and the output layer, as shown in Fig. \ref{fig:dnn}.
In the input layer, the neural network receives the three-point stencil $\vecS=(f_0, f_1, f_2)$ that is either $\left( f^+_{i-1}, f^+_i, f^+_{i+1} \right)$ in Fig. \ref{fig:stencil_l} or $\left( f^-_{i+2}, f^-_{i+1}, f^-_i \right)$ in Fig. \ref{fig:stencil_r}.
The Delta layer \eqref{eq:delta_layer}, which is used to pre-process the input stencil, is modified to
\begin{equation} \label{eq:delta_layer_modified}
\begin{aligned}
 & \tilde{\Delta}^\mdf_1 = \max(\left| f_0-f_1 \right|,\varepsilon), \, 
   \tilde{\Delta}^\mdf_2 = \max(\left| f_1-f_2 \right|,\varepsilon), \, 
   \tilde{\Delta}^\mdf_3 = \tilde{\Delta}_3, \, 
   \tilde{\Delta}^\mdf_4 = \tilde{\Delta}_4, \\
 & \Delta^\mdf_j = \tilde{\Delta}^\mdf_j / \max \left( \tilde{\Delta}^\mdf_1, \tilde{\Delta}^\mdf_2 \right), \,  \varepsilon^\mdf = 10^{-10}, \, j=1,2,3,4.
\end{aligned}     
\end{equation}
The modified Delta layer \eqref{eq:delta_layer_modified} exhibits enhanced stability compared to the original one \eqref{eq:delta_layer} in smooth regions.
To see this, consider the constant stencil $\vecS_1=(1,1,1)$ and the linear stencil $\vecS_2=(1,2,3)$.
Then the Delta layer yields
$$ \Delta(\vecS_1)=(0,0,0,0), \, \Delta(\vecS_2)=(1,1,2,0), $$
whereas for the modified Delta layer, we have
$$ \Delta^\mdf(\vecS_1)=(1,1,0,0), \, \Delta^\mdf(\vecS_2)=(1,1,2,0), $$
in which case, the modified Delta layer stably assigns the value of $1$ to the first two entries such that the model is better equipped to detect the smooth stencils.
Additionally, the modified Delta layer is scale-invariant if the stencil $\vecS=(f_0, f_1, f_2)$ and the scale factor $\kappa$ satisfy
$$
   \min \left\{ |f_0-f_1|, |f_1-f_2| \right\} \ge \varepsilon^\mdf, \, |\kappa| \min \left\{ |f_0-f_1|, |f_1-f_2| \right\} \ge \varepsilon^\mdf.
$$
However, the translation-invariance of the modified Delta layer holds regardless of the regularity of the stencil. 
After the modified Delta layer, the four features are processed by two hidden layers. 
Following \cite{Park}, each hidden layer consists of 16 nodes, and utilizes the GELU activation function.
The output layer is composed of 2 nodes and the softmax function as the activation function, which generates the nonlinear weights $\omega = (\omega_0, \omega_1) = \weno^\CADNN (\vecS)$.
The softmax function guarantees the convex combination of the numerical fluxes in the reconstruction procedure. 

\subsection{Dataset} \label{sec:fnn_dataset}
We now present the conservative approximation to the derivative using point values \cite{Shu}, which underpins the construction of the labeled dataset.

\subsubsection{Conservative approximation to derivative} \label{sec:fnn_dataset_cad}
Given the values of a function $v(x)$ at uniformly spaced points, $v_i = v(x_i), \, i=1,\cdots,N$, we want to find a numerical flux function $\hat{v}_{i+1/2} = \hat{v}(v_{i-1}, v_i, v_{i+1})$ for $i=2,\cdots,N-1$, such that the derivative $v'(x)$ is approximated by the flux difference, 
\begin{equation} \label{eq:conservative_approx}
 v'(x_i) \approx \frac{\hat{v}_{i+1/2}-\hat{v}_{i-1/2}}{\Delta x}, \, i=3,\cdots,N-1,
\end{equation}
where $\Delta x$ is the uniform spacing.
We wish to achieve third-order accuracy in smooth regions and avoid spurious oscillations when there is a discontinuity.
\begin{figure}[htbp]
\centering
\begin{tikzpicture}
\draw (-3,0) -- (3,0);
\filldraw[red] (-3,0) circle (2pt);
\filldraw[red] (-1,0) circle (2pt);
\filldraw[red] (1,0) circle (2pt);
\filldraw[red] (3,0) circle (2pt);
\node (11) [inner sep=1pt, minimum size=1pt] at (0,0.9)  {\hspace{0.3mm} \small $\hat{v}_{i-1/2}$};
\node (31) [inner sep=1.5pt, minimum size=1.5pt] at (1,1.3) {\hspace{-3mm} $v'(x_i)$};
\node (21) [inner sep=1pt, minimum size=1pt] at (2,0.9)  {\hspace{0.3mm} \small $\hat{v}_{i+1/2}$};
\node [rectangle, inner sep=2pt, minimum size=2pt, draw=green] (12) at (0,0) {};
\node [circle, inner sep=1.5pt, minimum size=1.5pt, draw=red] (32) at (1,0) {};
\node [rectangle, inner sep=2pt, minimum size=2pt, draw=blue] (22) at (2,0) {};
\draw [green,->,thick] (11)--(12);
\draw [blue,->,thick] (21)--(22);
\draw [black,->,thick] (31)--(32);
\node at (-3+0.12,-1/3) {$v_{i-2}$};
\node at (-1+0.12,-1/3) {$v_{i-1}$};
\node at (1,-1/3) {$v_{i}$};
\node at (3+0.12,-1/3) {$v_{i+1}$};
\draw (-3,-3/4-1/8+1/12) -- (1,-3/4-1/8+1/12);
\node at (-2,-2/4-1/8) {\footnotesize $\omega_{0,i-1/2}$};
\node at (0,-2/4-1/8) {\footnotesize $\omega_{1,i-1/2}$};
\filldraw[green] (-3,-3/4-1/8+1/12) circle (2pt);
\filldraw[green] (-1,-3/4-1/8+1/12) circle (2pt);
\filldraw[green] (1,-3/4-1/8+1/12) circle (2pt);
\draw (-1,-5/3+2/8+1/12) -- (3,-5/3+2/8+1/12);
\node at (0,-5/3+4/8) {\footnotesize $\omega_{0,i+1/2}$};
\node at (2,-5/3+4/8) {\footnotesize $\omega_{1,i+1/2}$};
\filldraw[blue] (-1,-5/3+2/8+1/12) circle (2pt);
\filldraw[blue] (1,-5/3+2/8+1/12) circle (2pt);
\filldraw[blue] (3,-5/3+2/8+1/12) circle (2pt);
\end{tikzpicture}
\caption{Conservative approximation to the derivative $v'(x_i)$ from a four-point stencil $(v_{i-2}, v_{i-1}, v_i, v_{i+1})$. The derivative $v'(x_{i-1})$ can also be approximated from the mirror-symmetric stencil $(v_{i+1}, v_i, v_{i-1}, v_{i-2})$.}
\label{fig:conservative_approx}
\end{figure}
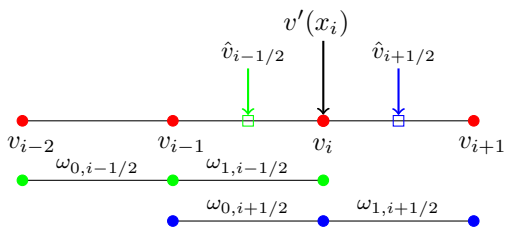
As in Section \ref{sec:weno}, the numerical flux function $\hat{v}_{i+1/2}$, defined on the stencil $(v_{i-1}, v_i, v_{i+1})$, takes a convex combination of two candidate numerical fluxes $\hat{v}^k_{i+1/2}$,
\begin{equation} \label{eq:numerical_flux_function}
\begin{gathered}
 \hat{v}_{i+1/2} = \omega_0 \hat{v}^0_{i+1/2} + \omega_1 \hat{v}^1_{i+1/2}, \\
 \hat{v}^0_{i+1/2} = -\frac{1}{2} v_{i-1} + \frac{3}{2} v_i, \quad \hat{v}^1_{i+1/2} = \frac{1}{2} v_i + \frac{1}{2} v_{i+1}.
\end{gathered}
\end{equation}
Consider the four-point stencil $(v_{i-2}, v_{i-1}, v_i, v_{i+1})$ in Fig. \ref{fig:conservative_approx}.
If the stencil is smooth and the weights $\omega_{k,i\pm1/2}$ coincide with the linear weights $d_k$, the flux difference approximation is third-order accurate,
\begin{equation} \label{eq:conservative_approx_third}
 \frac{\hat{v}_{i+1/2}-\hat{v}_{i-1/2}}{\Delta x} = v'(x_i) + O(\Delta x^3).   
\end{equation}
However, in the presence of a discontinuity within the stencil, the weight(s) corresponding to the discontinuous substencil(s) would be essentially zero, in accordance with the WENO methodology.
Similarly, the derivative $v'(x_{i-1})$ can be approximated from the mirror-symmetric stencil $(v_{i+1}, v_i, v_{i-1}, v_{i-2})$.

\subsubsection{Construction of input data and labels} \label{sec:fnn_dataset_datalabel} 
Based on the conservative approximation to the derivative, the dataset is constructed as follows.
Each sample in the dataset contains the four-point stencil and its label that is determined by the smoothness of the stencil.
Following the functions in \cite{Bezgin22}, we make use of a similar set of functions in Table \ref{tab:dataset}.
To generate data from smooth regions, we select (i) cubic polynomials with random coefficients, and (ii) sine functions and hyperbolic tangent functions with varying frequencies.
We extract the four-point stencils from those functions, using the discretized domain $[-1,1]$ with the uniform spacing $\Delta x=0.01$.
The label is defined as the exact derivative at the third point as shown in Fig. \ref{fig:conservative_approx}, as well as the derivative at the second point for the mirror-symmetric stencil.
We construct discontinuous four-point stencils from piecewise constant functions and linear polynomials with one jump discontinuity. 
Only one stencil per function, containing the discontinuity, is used in the dataset.
The label is defined as the average rate of change between two function values adjacent to the discontinuity.
\begin{table}[htbp]
\renewcommand{\arraystretch}{1.1}
\scriptsize
\centering
\caption{Functions for stencils in the dataset}      
\begin{tabular}{ccccc}
\hline 
Functions & Parameters & Number of samples \\
\hline         
$\sum_{i=0}^3 {\mathfrak a}_i x^i$                  & ${\mathfrak a}_i \sim \mathcal{U}[-1,1]$                    & 3920 \\
$\tanh({\mathfrak b} x), \sin({\mathfrak b} \pi x)$ & ${\mathfrak b} \sim \mathcal{U}[2,20]$                      & 7880 \\
${\mathfrak c}_0(x \le 0) + {\mathfrak c}_1(x>0)$   & ${\mathfrak c}_0, {\mathfrak c}_1 \sim \mathcal{U}[-10,10]$ & 8000 \\
$\pm x+{\mathfrak d}(x>0)$                          & ${\mathfrak d} \sim \mathcal{U}[0.5,2.5]$                   & 4000  \\
\hline
\end{tabular}
\label{tab:dataset}
\end{table}

\subsection{Training} \label{sec:fnn_training}
We employ the supervised learning to train the neural network using the Adam optimizer with the learning rate $10^{-4}$ and the weight decay $0.01$.

\subsubsection{Training pipeline} \label{sec:fnn_training_pipeline}
Consider one sample, the four-point stencil $(v_{i-2}, v_{i-1}, v_i, v_{i+1})$ and its label $Dv_i$, from the dataset.  
To approximate the derivative in the conservative form, the stencil is decomposed into two overlapping three-point substencils: $(v_{i-2}, v_{i-1}, v_i)$ and $(v_{i-1}, v_i, v_{i+1})$, as illustrated in Fig. \ref{fig:conservative_approx}.
Each substencil is provided as input to the feedforward neural network, which returns two sets of nonlinear weights: $(\omega_{0,i-1/2}, \omega_{1,i-1/2})$ and $(\omega_{0,i+1/2}, \omega_{1,i+1/2})$.  
Using these weights and the candidate fluxes in \eqref{eq:numerical_flux_function}, we compute the numerical fluxes $\hat{v}_{i-1/2}$ and $\hat{v}_{i+1/2}$ at the respective interfaces.
The conservative approximation to the derivative is then obtained via the flux difference:
$$
   D \hat{v}_i = \frac{\hat{v}_{i+1/2} - \hat{v}_{i-1/2}}{\Delta x},
$$
which serves as the prediction.
During the training phase, the neural network is optimized to match the prediction with the label $Dv_i$.

\subsubsection{Loss function} \label{sec:fnn_training_loss}
To design the feedforward neural network for WENO weighting function, we consider the following key properties: 
\begin{enumerate}[label=\arabic*., font=\itshape]
\item For smooth stencils, the nonlinear weights $\omega_k$ are expected to recover the linear weights $d_k$, ensuring the spatial third-order accuracy \eqref{eq:conservative_approx_third}.
\item For the stencil containing a discontinuity, it is hoped that the nonlinear weight $\omega_k$ associated with the discontinuous substencil is essentially zero for the ENO behavior. 
\item Given the set of nonlinear weights $\omega=(\omega_0,\omega_1)$ for the stencil $\vecS = (f_0,f_1,f_2)$, the flipped stencil $\breve{\vecS}=(f_2,f_1,f_0)$ corresponds to the set of nonlinear weights $\breve{\omega}$, 
\begin{equation} \label{eq:symmetry}
 \breve{\omega} = (\breve{\omega}_0,\breve{\omega}_1) = \left( \frac{\omega_1}{4\omega_0+\omega_1}, \frac{4\omega_0}{4\omega_0+\omega_1} \right).
\end{equation}
\end{enumerate}
The first two properties are crucial to emulate the WENO idea, while WENO3-JS and WENO3-Z weighting functions share the third symmetric property.
Accordingly, we define the loss function as
\begin{equation} \label{eq:loss}
 \loss = \loss_\CAD + C \cdot \loss_\SYM + D \cdot \loss_\LN,
\end{equation}
where $C$ and $D$ are tunable hyperparameters.
The loss function comprises three components: the loss of conservative approximation to derivative $\loss_\CAD$, the symmetric-balancing term $\loss_\SYM$ matching the third property, and the linear term $\loss_\LN$. 

\begin{figure}[htbp]
\centering
\includegraphics[width=0.325\textwidth]{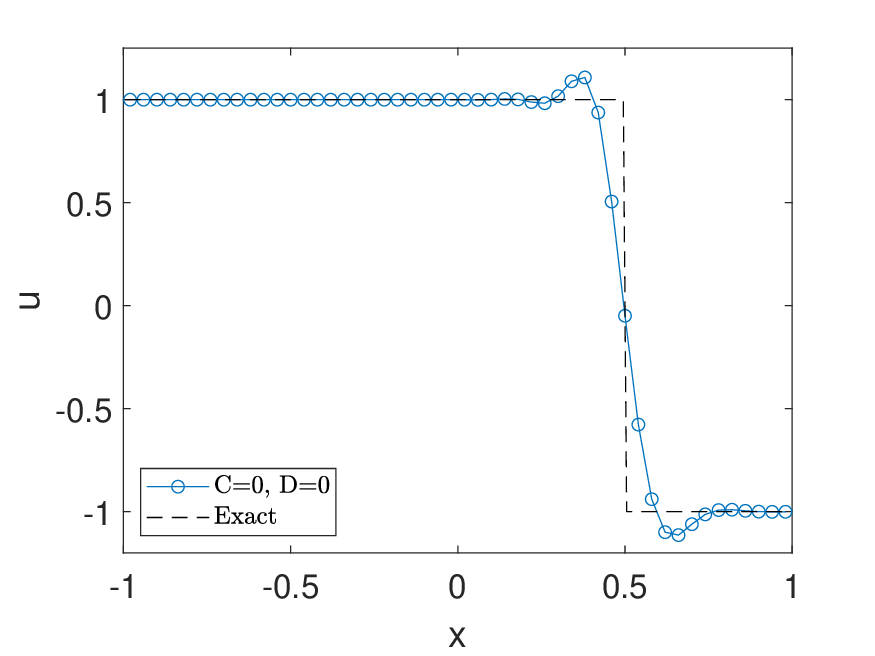}
\includegraphics[width=0.325\textwidth]{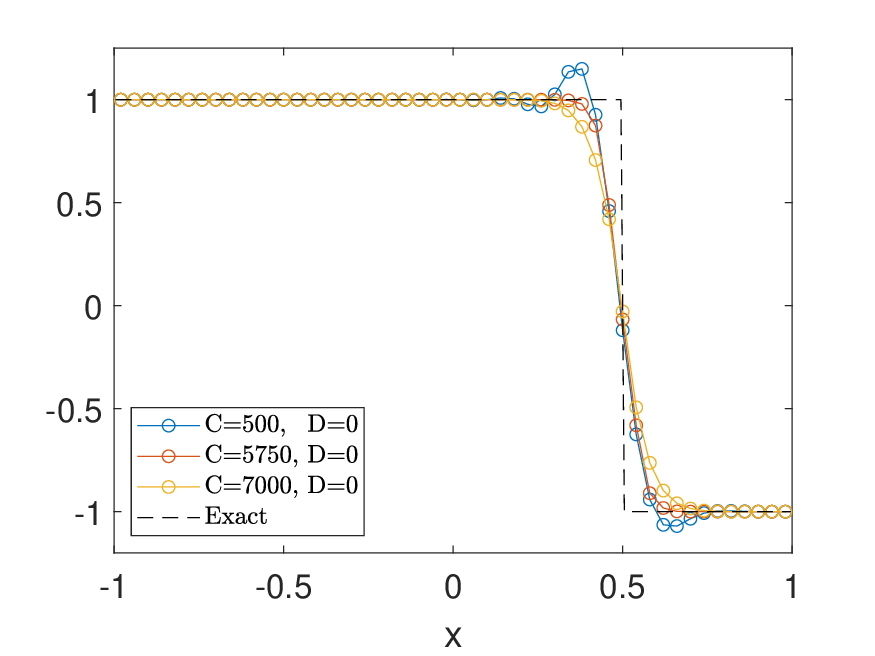}
\includegraphics[width=0.325\textwidth]{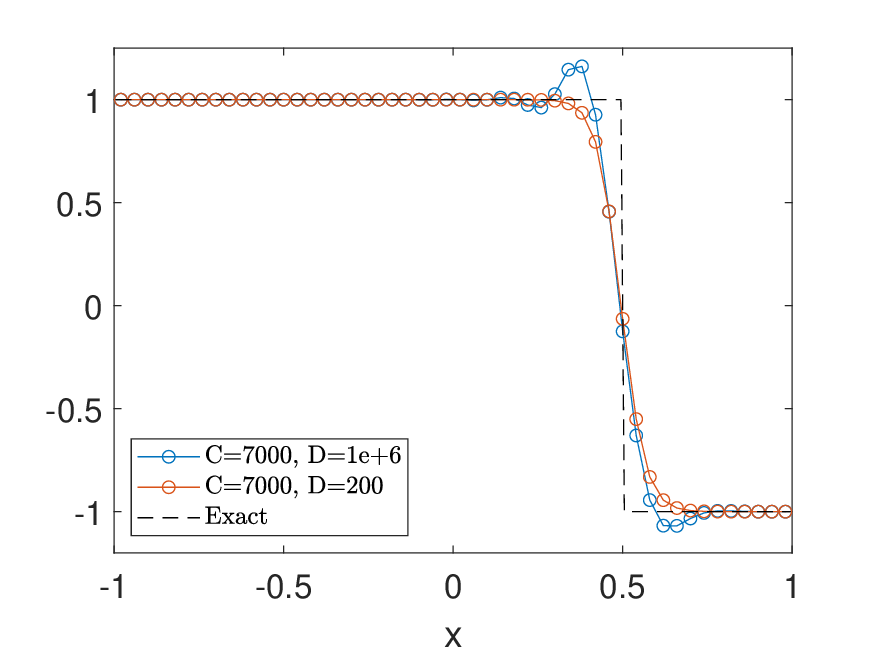}
\caption{Solution profiles at $T=0.5$ for the Riemann problem of one-dimensional advection equation \eqref{eq:advection_1d}, computed by the feedforward neural network with the loss function \eqref{eq:loss}: $C=0$ and $D=0$ (left), $C=500, 5750, 7000$ and $D=0$ (middle), and $C=7000$ and $D=200, 10^6$ (right).}
\label{fig:hyperparameter}
\end{figure}

The first term $\loss_\CAD$, 
$$ \loss_\CAD = \frac{1}{N_b} \sum_{l=1}^{N_b} \left( D \hat{v}_l-Dv_l \right)^2, $$
measures the squared difference between the conservative approximation to the derivative $D \hat{v}_l$ using the predicted nonlinear weights from the neural network and the the labeled reference $Dv_l$, over a mini-batch of size $N_b=200$.
However, as shown in the left panel of Fig. \ref{fig:hyperparameter}, using only this loss term $\loss_\CAD$ ($C=D=0$) results in oscillatory solutions near the discontinuity in the Riemann problem for the linear advection equation \eqref{eq:advection_1d} at the final time $T=0.5$.
This motivates the inclusion of some additional regularizing term to suppress oscillations.

For each mini-batch of $N_b$ four-point stencils, we have totally $2N_b$ three-point stencils. 
For the $l$th stencil $\vecS_l=(f_{0,l},f_{1,l},f_{2,l})$, the neural network gives the nonlinear weights, 
$$ \omega^\CADNN_l = ( \omega^\CADNN_{0,l}, \omega^\CADNN_{1,l} ) = \weno^\CADNN (\vecS_l). $$
The flipped stencil $\breve{\vecS}_l=(f_{2,l},f_{1,l},f_{0,l})$ yields
$$ \breve{\omega}^\CADNN_l = ( \breve{\omega}^\CADNN_{0,l}, \breve{\omega}^\CADNN_{1,l} ) = \weno^\CADNN (\breve{\vecS}_l). $$ 
According to the third property, the expected nonlinear weights for the flipped stencil $\breve{\vecS}_l$ are
$$ 
   \check{\omega}_l = \left( \frac{\omega^\CADNN_{1,l}}{4\omega^\CADNN_{0,l}+\omega^\CADNN_{1,l}}, \frac{4\omega^\CADNN_{0,l}}{4\omega^\CADNN_{0,l}+\omega^\CADNN_{1,l}} \right). 
$$
The symmetric-balancing term $\loss_\SYM$ is then defined as
$$ \loss_\SYM = \frac{1}{N_b} \sum_{l=1}^{2N_b} \left\| \log{\breve{\omega}^\CADNN_l}-\log{\check{\omega}_l} \right\|^2_2. $$
We use the mean squared logarithmic error in \cite{Park} for better performance.
The hyperparameter $C$ in \eqref{eq:loss} quantifies how much the oscillations are suppressed.
As shown in the middle panel of Fig. \ref{fig:hyperparameter}, increasing $C$ enforces the symmetric constraint and suppresses oscillations near discontinuities. However, excessive emphasis on symmetry may introduce dissipation.

To further reduce the numerical dissipation around discontinuities while retaining third-order accuracy in smooth regions, we add the linear term $\loss_\LN$ to the loss function, as we did in \cite{Park}.
We begin by modifying the smoothness gauge \cite{Park} for the three-point stencil, which is defined by
$$
   \lambda_l = \e^{-6 r_l}, \quad r_l = \max \left( \frac{\Delta^\mdf_{1,l}}{\Delta^\mdf_{2,l}}, \frac{\Delta^\mdf_{2,l}}{\Delta^\mdf_{1,l}} \right),
$$
where $\Delta^\mdf_{1,l}$ and $\Delta^\mdf_{2,l}$ are the first two entries of the modified Delta layer $\Delta^\mdf(S_l)$ \eqref{eq:delta_layer_modified}.
When $S_l$ is smooth, the ratio $r_l \approx 1$, and consequently, $\lambda_l$ approaches the positive constant $\e^{-6}$.
In contrast, for $S_l$ containing a discontinuity, $r_l$ becomes significantly greater than $1$, causing $\lambda_l$ to decay toward zero.
We adopt the definition of the linear term $\loss_\LN$ in \cite{Park},
$$
     \loss_\LN = \frac{1}{N_b} \sum_{l=1}^{2N_b} \lambda_l \left[ \log \left( 2 \omega^\CADNN_{0,l} \right) - \log \left( \omega^\CADNN_{1,l} \right) \right]^2,
$$
with $2N_b$ three-point stencils out of $N_b$ four-point stencils.
For the smooth stencil, the smoothness gauge $\lambda_l$ remains non-negligible, emphasizing the linear term $\loss_\LN$ in the optimization objective. 
This encourages the output weights to align closely with the linear weights.
Conversely, when the stencil contains a discontinuity, the smoothness gauge $\lambda_l\approx 0$, effectively deactivates the influence of the linear term in the loss function.
The hyperparameter $D$ controls the strength of this linearity enforcement. 
The right panel of Fig. \ref{fig:hyperparameter} illustrates the performance of this feedforward neural network for two different values of the hyperparameter $D$, with $C$ fixed, in the simulation of the Riemann problem of the linear advection problem \eqref{eq:advection_1d} at $T=0.5$. 
We can see that moderate values of $D$ (e.g., $D=200$) effectively reduce dissipation without introducing oscillations, but overly large $D$ (e.g., $D = 10^6$) leads to oscillatory behavior around the discontinuity.

Figure \ref{fig:hyperparameter} shows the effect of varying $C$ and $D$ on the simulation quality for the Riemann problem of the linear advection problem.
The levels of numerical dissipation and oscillation can be finely controlled via the hyperparameters $C$ and $D$ in the loss function \eqref{eq:loss}. 
Based on the these results, we choose the two combinations of the hyperparameters $(C, D) = (5750, 0)$ and $(7000, 800)$, which are referred to as WENO3-CADNN1 and WENO3-CADNN2, respectively, for the numerical examples in Section \ref{sec:nr}.

\section{Numerical results} \label{sec:nr}
In this section, we carry out several one- and two-dimensional numerical experiments to show the performance of WENO3-CADNN1 and WENO3-CADNN2, in comparison with the classical WENO3-Z and WENO5-JS schemes. 
The explicit third-order TVD Runge-Kutta method \cite{ShuOsherI,Gottlieb} is used for time integration with a CFL number of $0.4$.

\subsection{One-dimensional scalar problem}
\begin{figure}[htbp]
\centering
\includegraphics[height=0.33\textwidth]{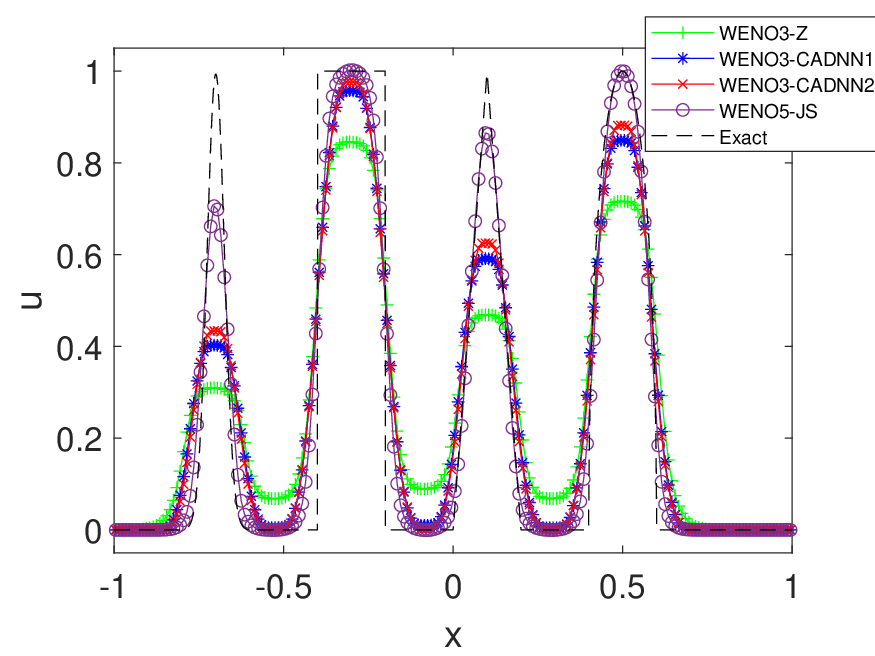}
\includegraphics[height=0.33\textwidth]{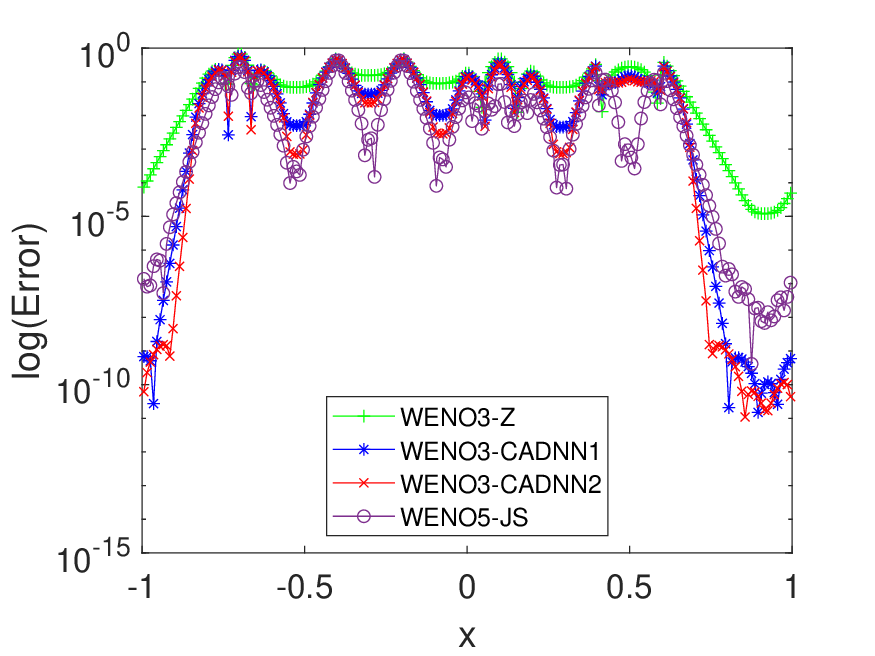}
\caption{Solution profiles (left) and log-scale pointwise errors (right) for Example \ref{ex:advection_1d} at $T=8$ approximated by WENO3-Z (green), WENO3-CADNN1 (blue), WENO3-CADNN2 (red) and WENO5-JS (purple) with $N = 200$. 
The dashed black line is the exact solution.}
\label{fig:advection_1d}
\end{figure}

\begin{example} \label{ex:advection_1d}
Consider the one-dimensional advection equation
\begin{equation} \label{eq:advection_1d}
 u_t + u_x = 0,
\end{equation} 
with the initial condition
\begin{align*}
 & u(x,0) = \left\{ 
   			 \begin{array}{ll} 
   			  \frac{1}{6} \left[ G(x, \beta, z-\delta) + 4 G(x, \beta, z) + G(x, \beta, z+\delta)\right], & -0.8 \leqslant x \leqslant -0.6, \\
			      1, & -0.4 \leqslant x \leqslant -0.2, \\
			      1 - \left| 10(x-0.1) \right|, & 0 \leqslant x \leqslant 0.2, \\
			      \frac{1}{6} \left[ F(x, \alpha, y-\delta) + 4 F(x, \alpha, y) + F(x, \alpha, y+\delta)\right], & 0.4 \leqslant x \leqslant 0.6, \\
			      0, & \text{otherwise},
			     \end{array} 
           \right. \\
 & G(x, \beta, z) = \e^{-\beta(x-z)^2}, \, F(x, \alpha, y) = \sqrt{\max \left( 1-\alpha^2 (x-y)^2, 0 \right)}, \\
 & \delta = 0.005, \, \beta = \frac{\ln 2}{36 \delta^2}, \, z = -0.7, \, \alpha = 10, \, y = 0.5.
 \end{align*}
The computational domain $[-1, 1]$ is discretized to $N=200$ uniform cells. 
We run the simulation up to the final time $T=8$, with the time step $\Delta t = \cfl \cdot \Delta x$.
Fig. \ref{fig:advection_1d} shows the numerical solutions and pointwise errors on a logarithmic scale at the cell centers, which illustrates that WENO3-CADNNs outperform WENO3-Z within the class of third-order WENO schemes, but exhibit lower accuracy than WENO5-JS due to their lower formal order of accuracy.
\end{example}

\subsection{One-dimensional system problems}
For one-dimensional system problems, we consider the one-dimensional compressible Euler equations
\begin{equation} \label{eq:euler_1d}
 \bfu_t + \bff(\bfu)_x = 0, 
\end{equation}
where the column vector $\bfu$ of conserved variables and the flux vector $\bff$ in the $x$ direction are
$$
   \bfu = \left[ \rho, \, \rho u, \, E \right]^T, \quad \bff(\bfu) = \left[ \rho u, \, \rho u^2 + P, \, u(E+P) \right]^T,
$$
with $\rho$, $u$ and $P$ corresponding to the density, velocity and pressure, respectively.
The specific kinetic energy $E$ is defined as
$$
   E = \frac{P}{\gamma - 1} + \frac{1}{2} \rho u^2,
$$
where $\gamma = 1.4$ is given for the ideal gas.

\begin{figure}[htbp]
\centering
\includegraphics[height=0.35\textwidth]{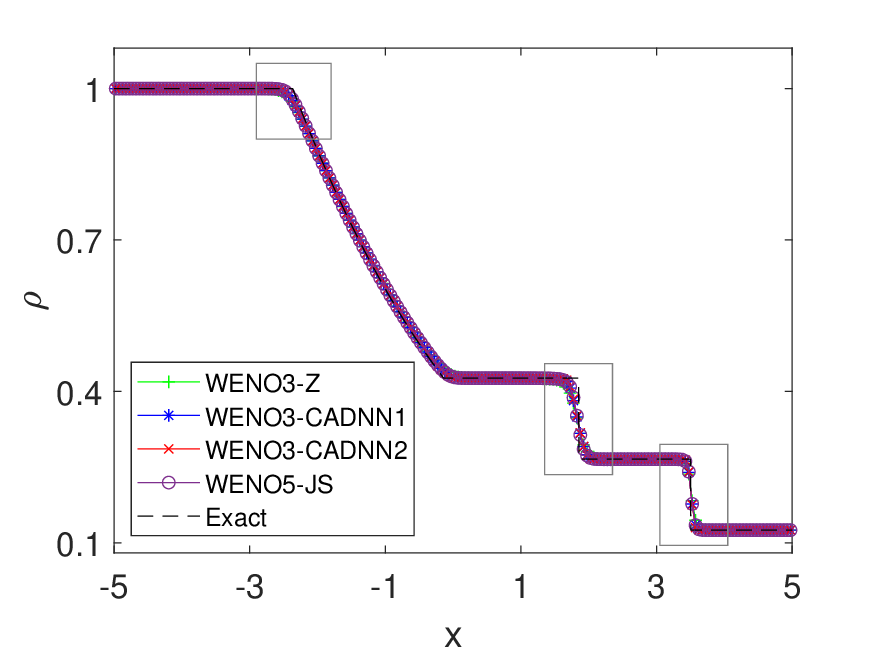} 
\includegraphics[height=0.35\textwidth]{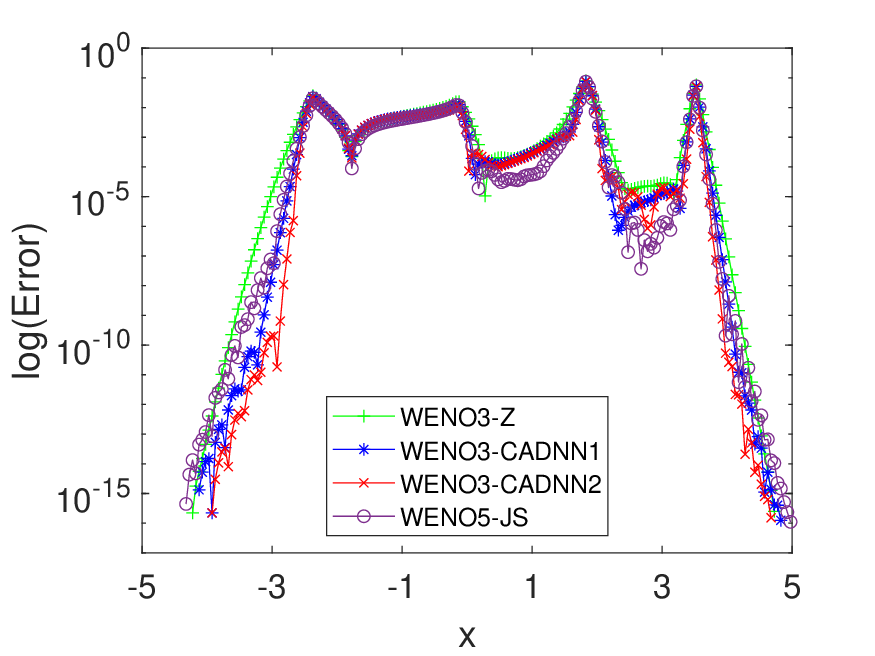} \\
\includegraphics[height=0.23\textwidth]{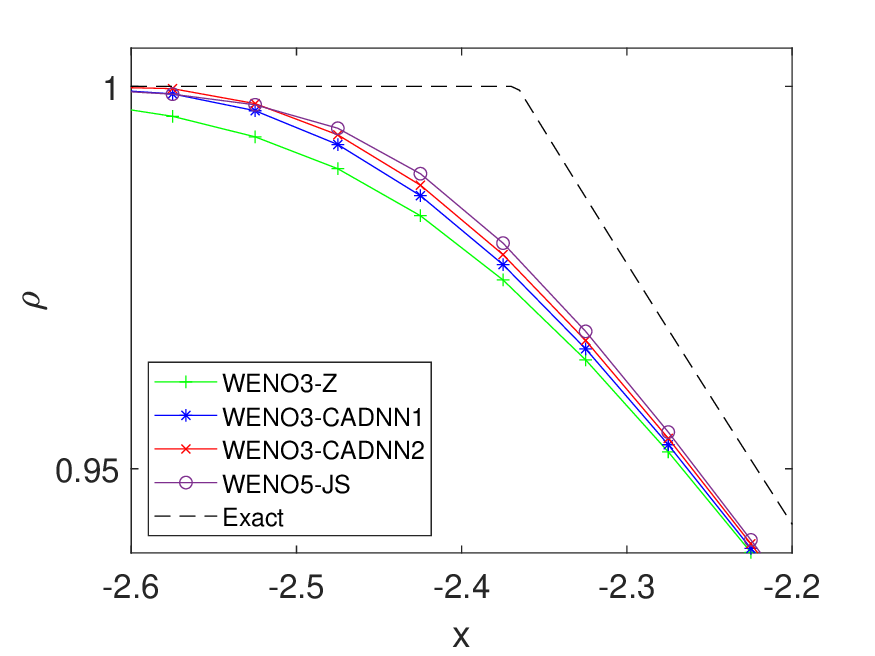} 
\includegraphics[height=0.23\textwidth]{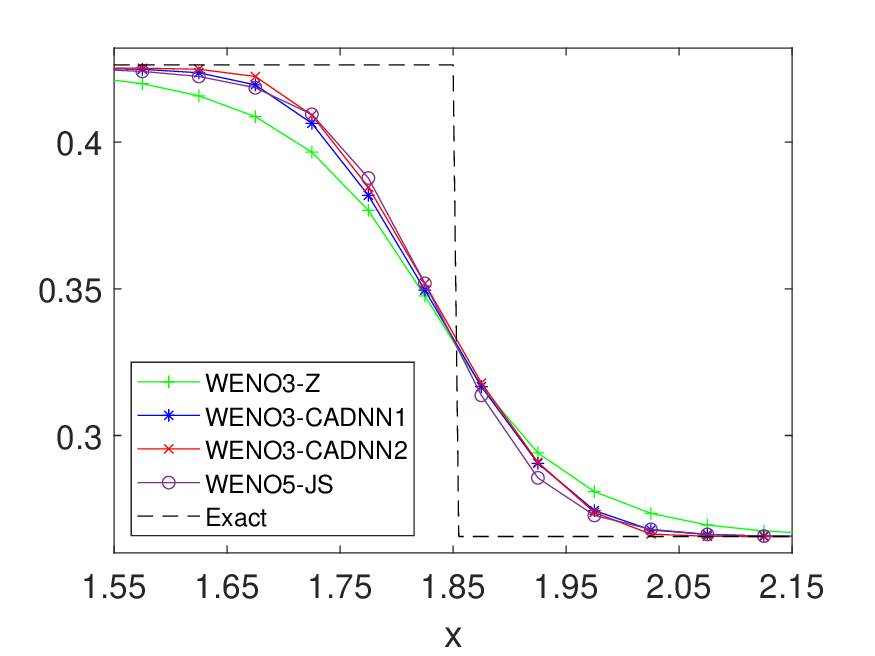}
\includegraphics[height=0.23\textwidth]{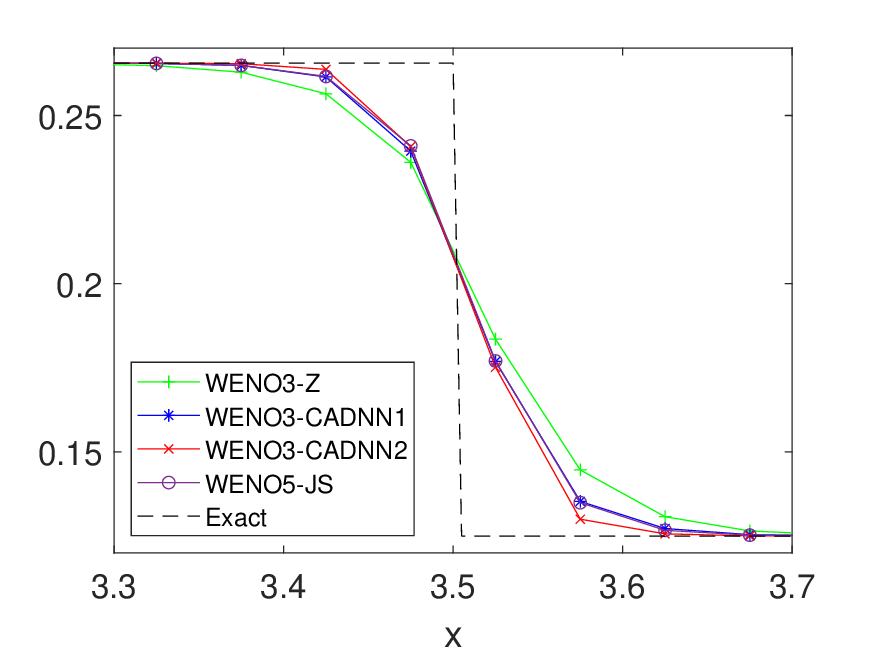}
\caption{Density profiles for the Sod problem \eqref{eq:euler_1d} and \eqref{eq:sod} at $T=2$ (top left), log-scale pointwise error (top right) and close-up view of the solutions in the boxes from left to right (bottom left, bottom middle, bottom right) approximated by WENO3-Z (green), WENO3-CADNN1 (blue), WENO3-CADNN2 (red) and WENO5-JS (purple) with $N = 200$. 
The dashed black line is the exact solution.}
\label{fig:sod}
\end{figure}

\begin{figure}[htbp]
\centering
\includegraphics[height=0.35\textwidth]{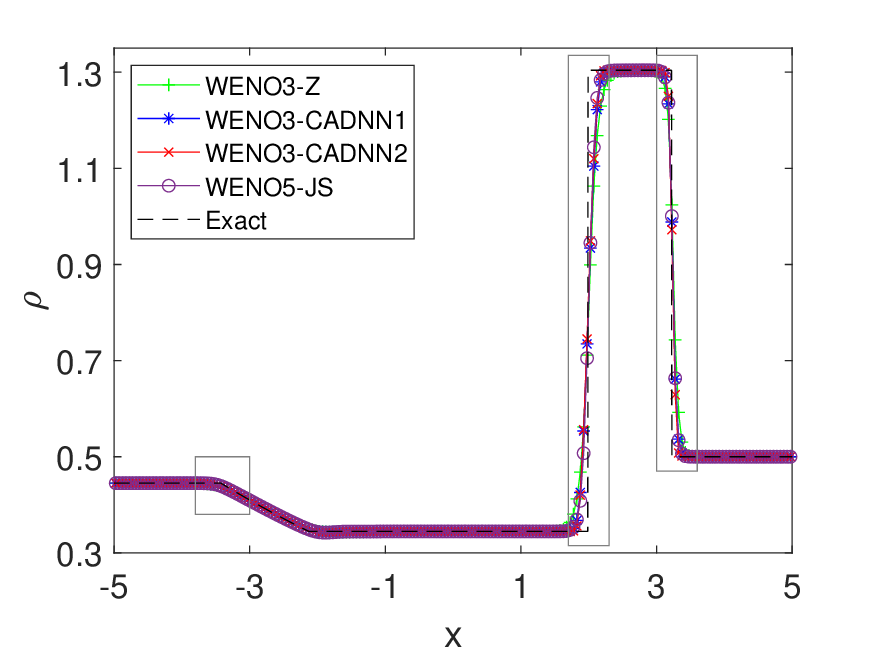}
\includegraphics[height=0.35\textwidth]{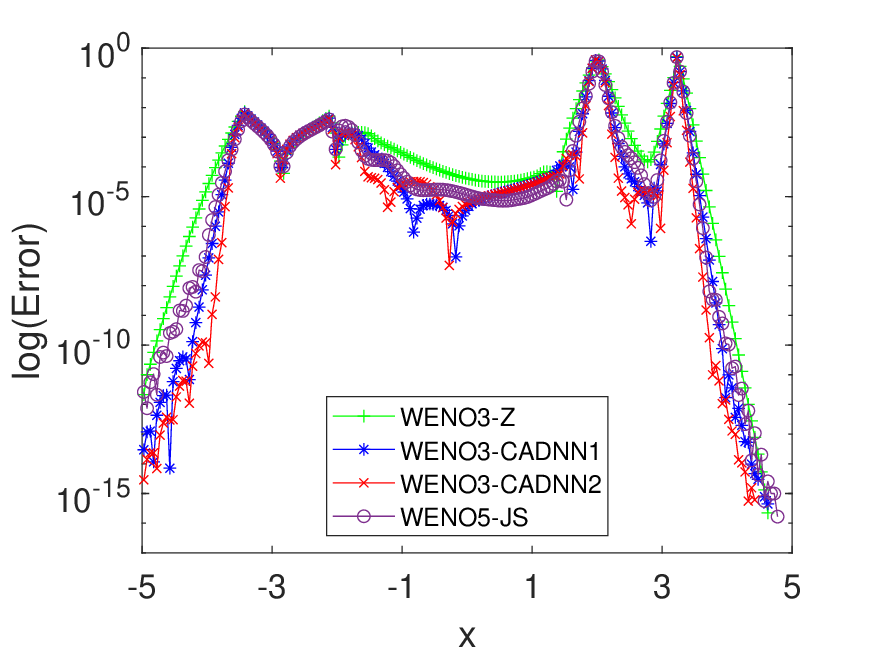}
\includegraphics[height=0.23\textwidth]{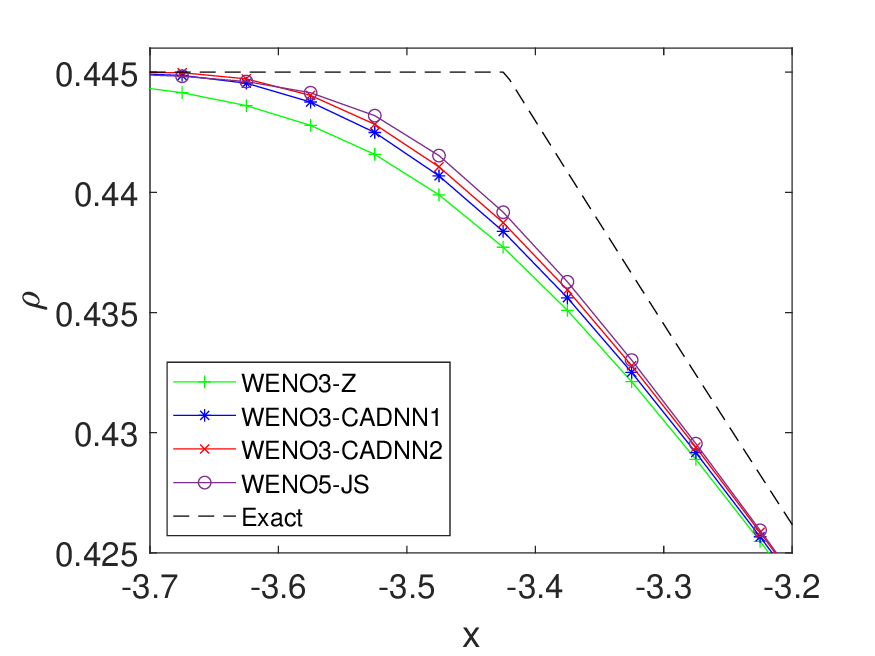}
\includegraphics[height=0.23\textwidth]{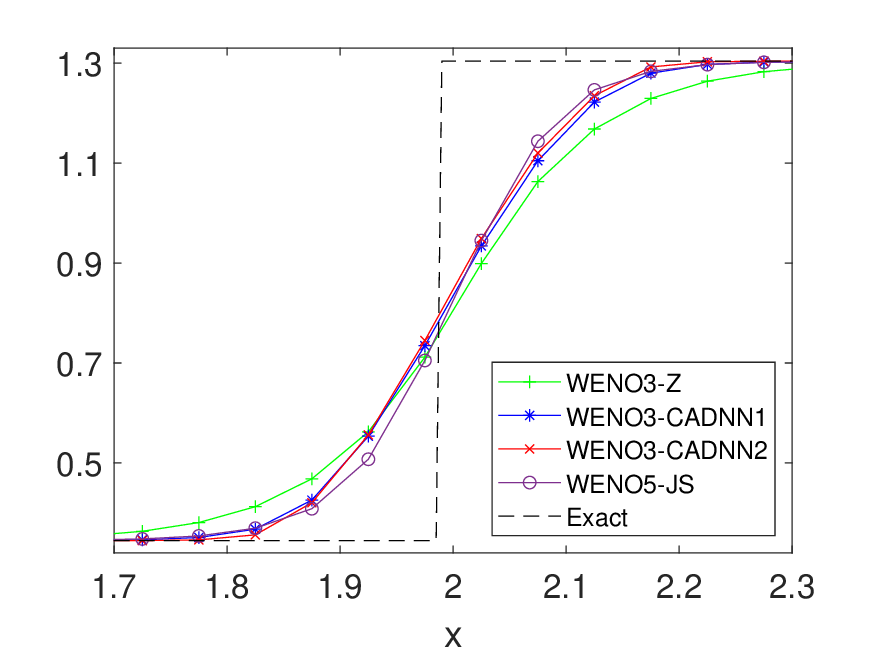}
\includegraphics[height=0.23\textwidth]{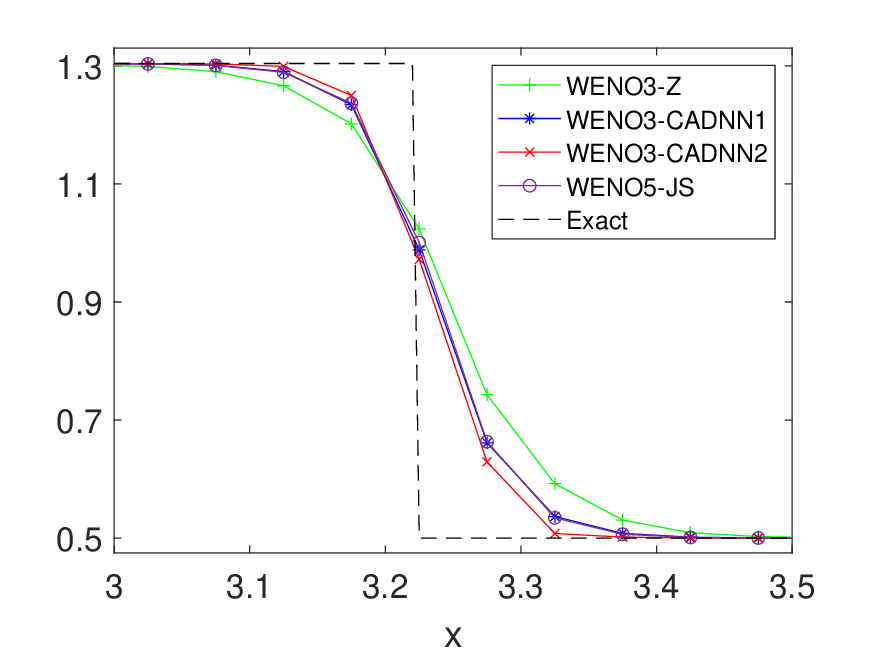}
\caption{Density profiles for the Lax problem \eqref{eq:euler_1d} and \eqref{eq:lax} at $T=1.3$ (top left), log-scale pointwise error (top right) and close-up view of the solutions in the boxes from left to right (bottom left, bottom middle, bottom right) approximated by WENO3-Z (green), WENO3-CADNN1 (blue), WENO3-CADNN2 (red) and WENO5-JS (purple) with $N = 200$.
The dashed black line is the exact solution.}
\label{fig:lax}
\end{figure}

\begin{figure}[htbp]
\centering
\includegraphics[height=0.35\textwidth]{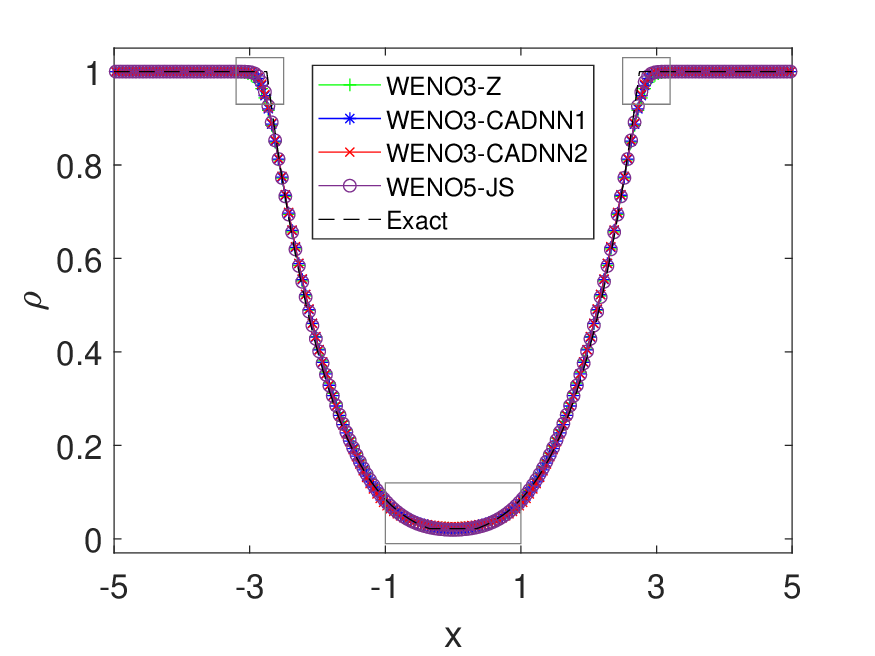}
\includegraphics[height=0.35\textwidth]{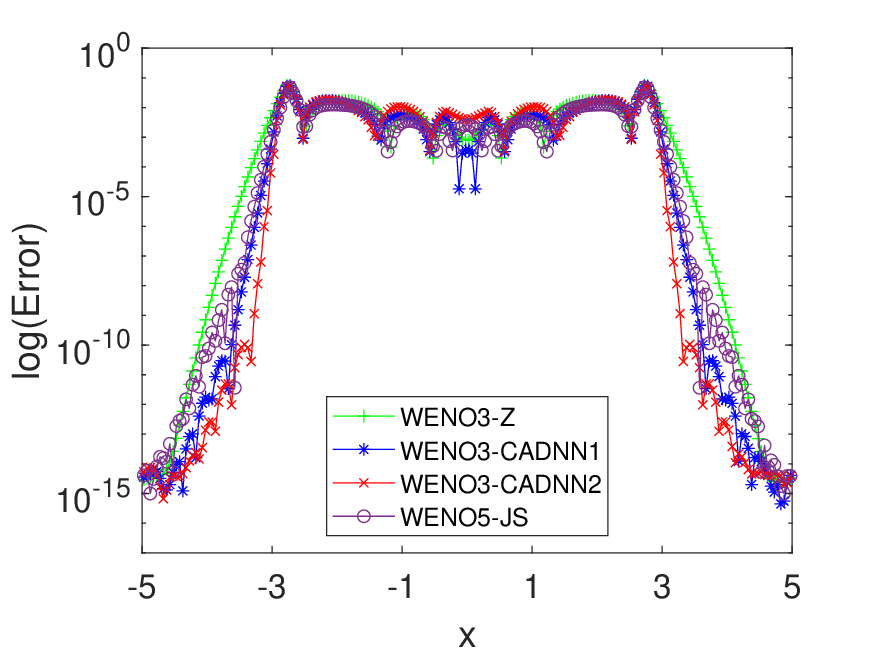}
\includegraphics[height=0.23\textwidth]{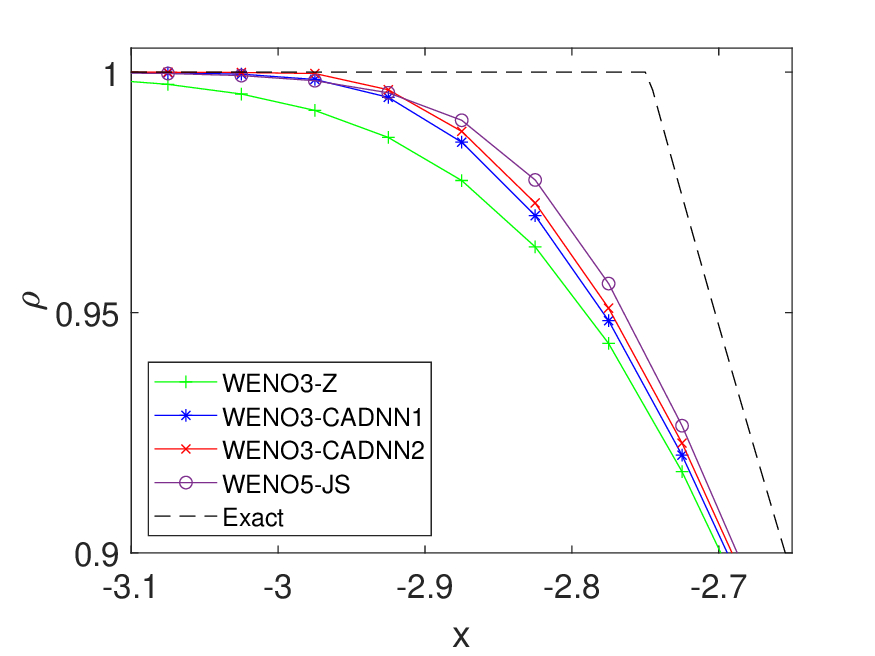}
\includegraphics[height=0.23\textwidth]{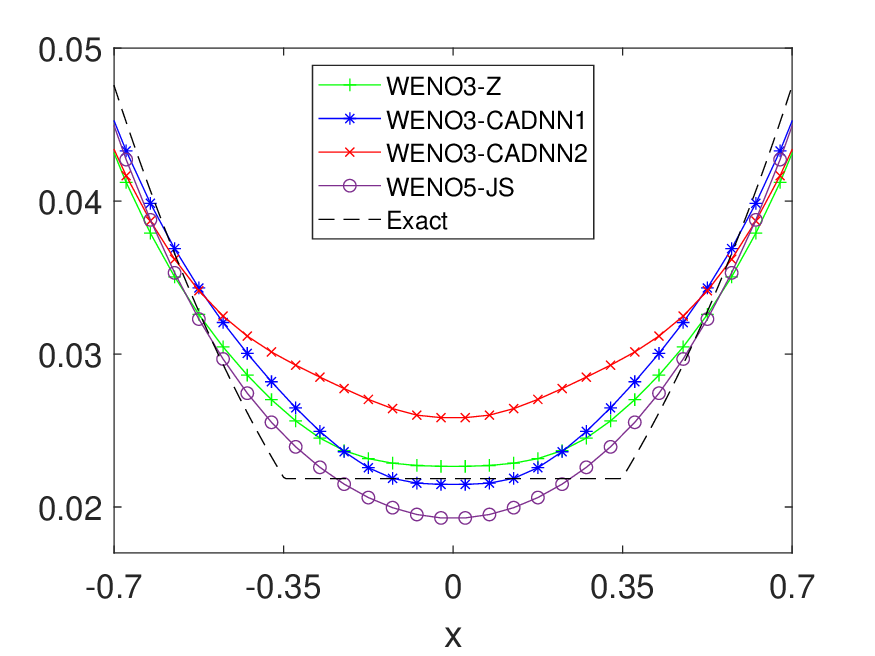}
\includegraphics[height=0.23\textwidth]{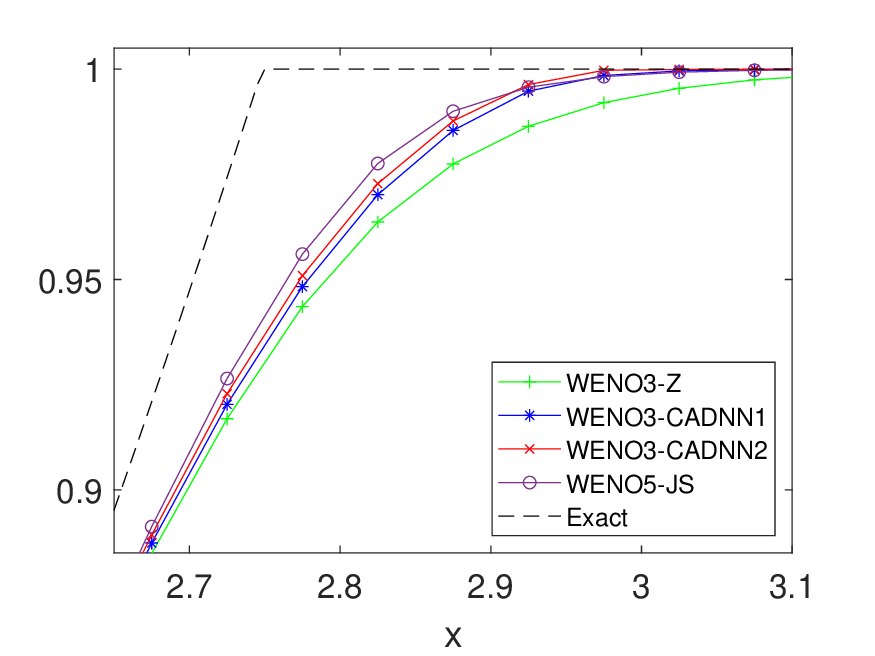}
\caption{Density profiles for the 123 problem \eqref{eq:euler_1d} and \eqref{eq:one23} at $T=1$ (top left), log-scale pointwise error (top right) and close-up view of the solutions in the boxes from left to right (bottom left, bottom middle, bottom right) approximated by WENO3-Z (green), WENO3-CADNN1 (blue), WENO3-CADNN2 (red) and WENO5-JS (purple) with $N = 200$. 
The dashed black line is the exact solution.}
\label{fig:one23}
\end{figure}

\begin{figure}[htbp]
\centering
\includegraphics[height=0.32\textwidth]{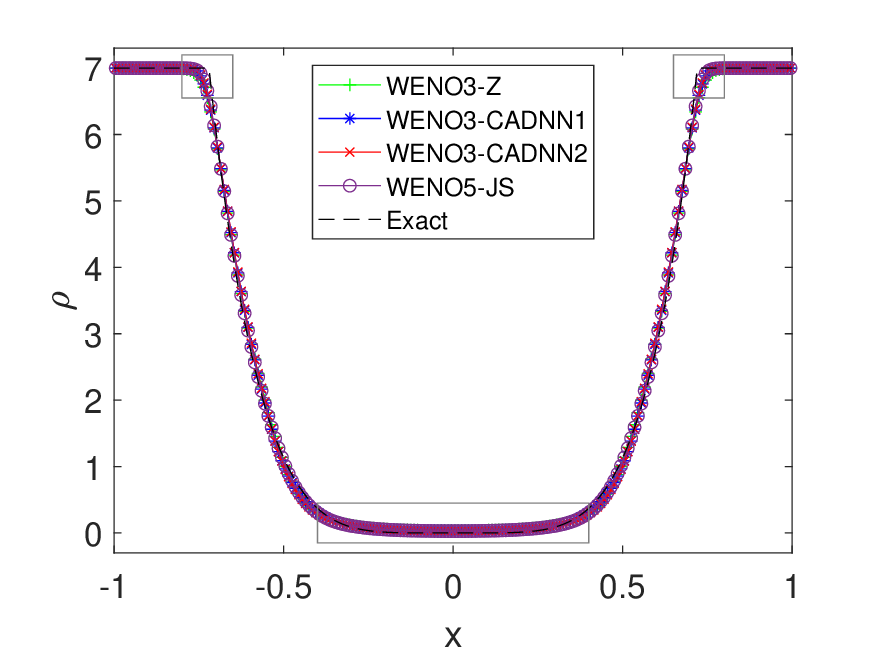} 
\includegraphics[height=0.32\textwidth]{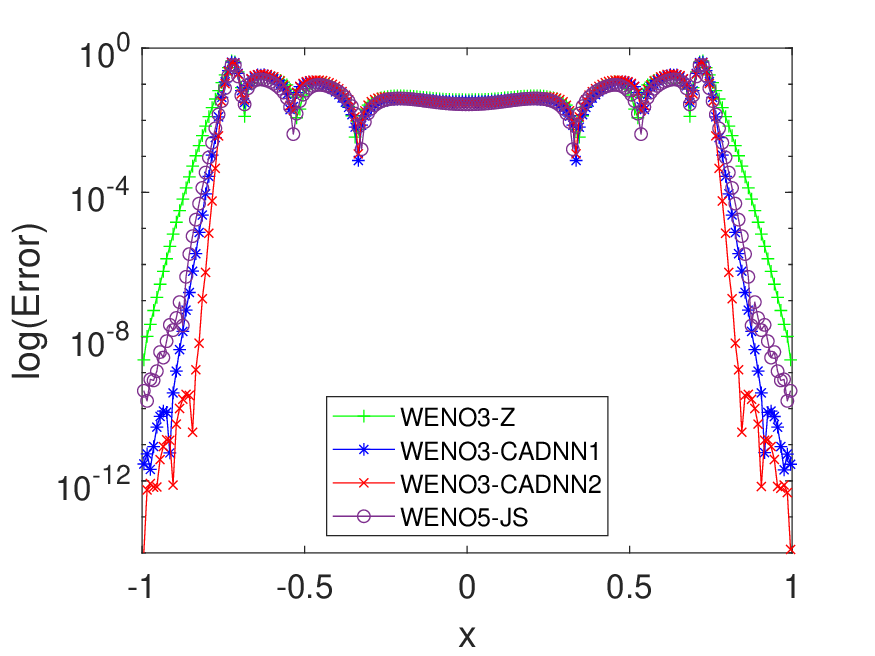}
\includegraphics[height=0.22\textwidth]{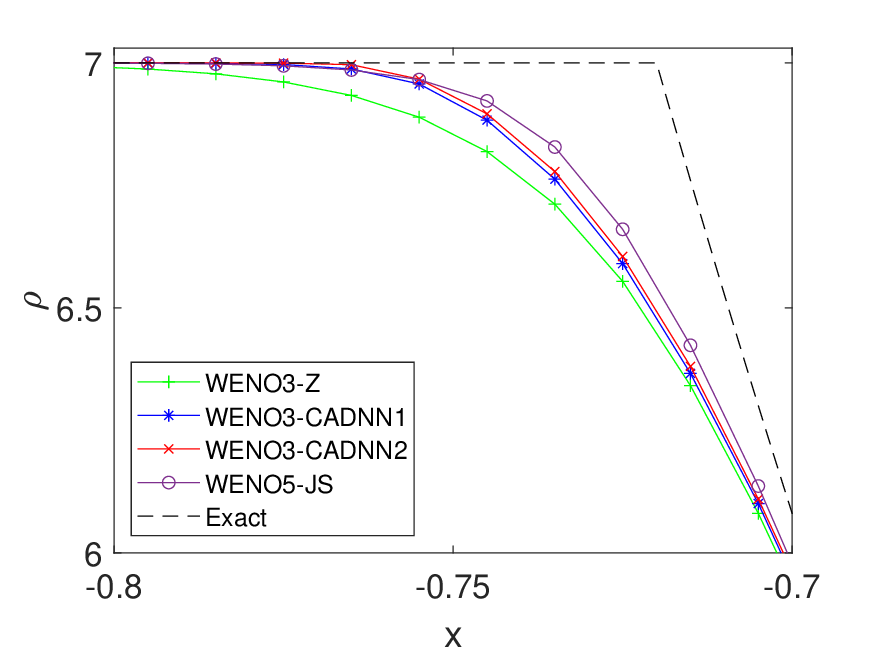} 
\includegraphics[height=0.22\textwidth]{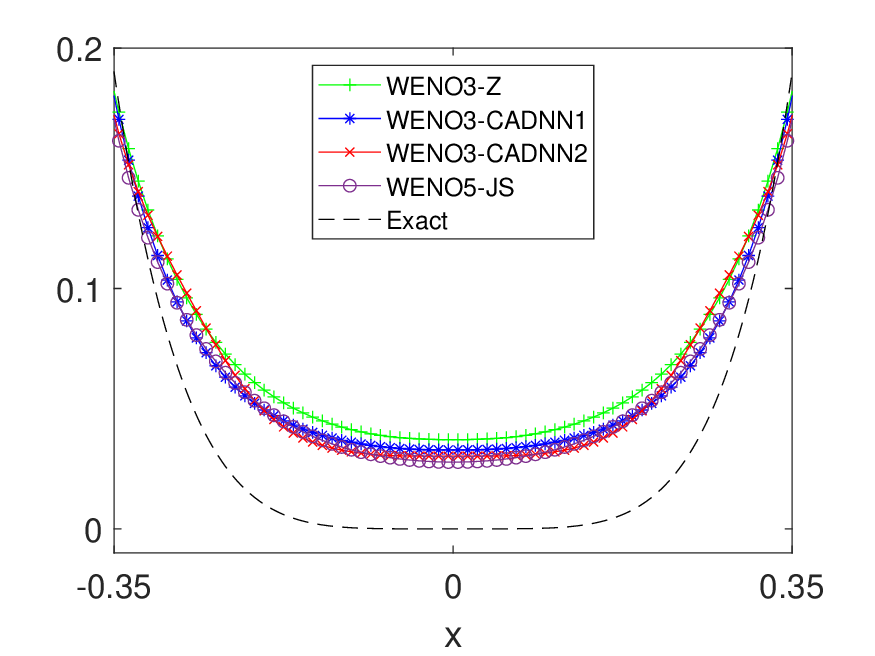}
\includegraphics[height=0.22\textwidth]{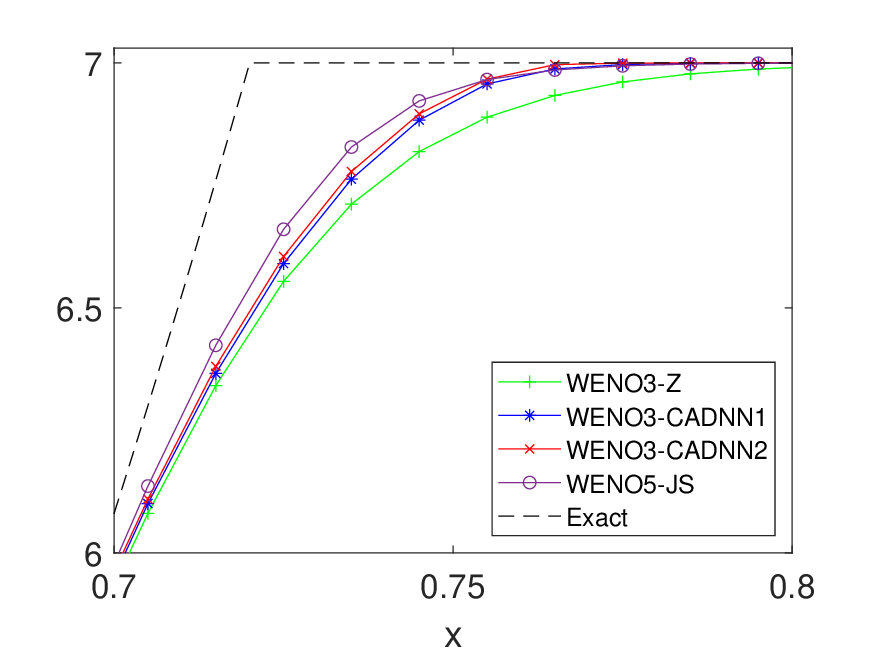}
\vspace{-0.3cm}
\caption{Density profiles for the double rarefaction problem \eqref{eq:euler_1d} and \eqref{eq:double_rarefaction} at $T=0.6$ (top left), log-scale pointwise error (top right) and close-up view of the solutions in the boxes from left to right (bottom left, bottom middle, bottom right) approximated by WENO3-Z (green), WENO3-CADNN1 (blue), WENO3-CADNN2 (red) and WENO5-JS (purple) with $N = 200$. 
The dashed black line is the exact solution.}
\label{fig:double_rarefaction}
\end{figure}

\begin{example} \label{ex:euler_1d}
We consider Riemann problems for the Euler equations, for which the exact solution is available to evaluate the accuracy of those WENO schemes.

As a benchmark shock tube test case, the Sod problem is \eqref{eq:euler_1d} with the initial condition
\begin{equation} \label{eq:sod}
 (\rho, u, P ) = \left\{ 
                  \begin{array}{ll} 
                   (1, \, 0, \, 1),       & x \leqslant 0, \\ 
                   (0.125, \, 0, \, 0.1), & x > 0.
                  \end{array} 
                 \right.
\end{equation}
The computational domain $[-5, 5]$ is discretized with $N=200$ uniform cells.
In Fig. \ref{fig:sod}, we show the density profile from each WENO scheme until $T=2$, along with the pointwise errors on a logarithmic scale.
The density profiles computed by WENO3-CADNNs exhibit reduced numerical dissipation around the regions of rarefaction, contact discontinuity and shock wave compared to WENO3-Z.
In particular, WENO3-CADNN2 yields more accurate results than WENO5-JS near the shock wave.

The Lax problem is another benchmark shock tube test case, for which the initial condition is 
\begin{equation} \label{eq:lax}
 (\rho, u, P ) = \left\{ 
                  \begin{array}{ll} 
                   (0.445, \, 0.698, \, 3.528), & x \leqslant 0, \\ 
                   (0.5, \, 0, \, 0.571),       & x > 0. 
                  \end{array} 
                 \right. 
\end{equation}
The computational domain is $[-5, 5]$ with $N=200$ uniform grid points. 
In Fig. \ref{fig:lax}, we plot the numerical solutions of the density at the final time $T=1.3$ and the corresponding log-scale pointwise errors at the grid points.
It is observed that WENO3-CADNNs show improved performance over WENO3-Z.
Among all, WENO3-CADNN2 provides the most accurate density profile around the region of shock wave.

The 123 problem also serves as a test case for the shock tube, with the initial condition
\begin{equation} \label{eq:one23}
 (\rho, u, P ) = \left\{ 
                  \begin{array}{ll} 
                   (1, \, -2, \, 0.4), & x \leqslant 0, \\ 
                   (1, \, 2, \, 0.4),  & x > 0. 
                  \end{array} 
                 \right. 
\end{equation}
We divide the computational domain $[-5, 5]$ into $N=200$ uniform cells.
Fig. \ref{fig:one23} shows the numerical density $\rho$ at $T=1$ with the pointwise errors on a logarithmic scale.
In the vicinity of the two moving rarefaction waves, WENO3-CADNNs provide more accurate density profiles than WENO3-Z, but are less accurate than WENO5-JS.
Near the stationary contact discontinuity, WENO3-CADNN1 yields the most accurate density approximation among the compared schemes.

The double rarefaction problem is characterized by the formation of a vacuum region in the center.
The initial condition is given by
\begin{equation} \label{eq:double_rarefaction}
 (\rho, u, P ) = \left\{
                  \begin{array}{ll} 
                   (7, \, -1, \, 0.2), & x \leqslant 0, \\ 
                   (7, \, 1, \, 0.2),  & x > 0,
                  \end{array} 
                 \right.
\end{equation}
on the computational domain $[-1, 1]$ with $N=200$ uniformly spaced grid points.
The numerical density profiles at the final time $T=0.6$, along with the pointwise errors plotted on a logarithmic scale, are presented in Fig. \ref{fig:double_rarefaction}.
Across all regions including the two rarefaction waves and the central vacuum zone, both WENO3-CADNN schemes exhibit improved density approximations compared to WENO3-Z, but offer lower resolution than WENO5-JS.
\end{example}

\begin{figure}[htbp]
\centering
\includegraphics[height=0.32\textwidth]{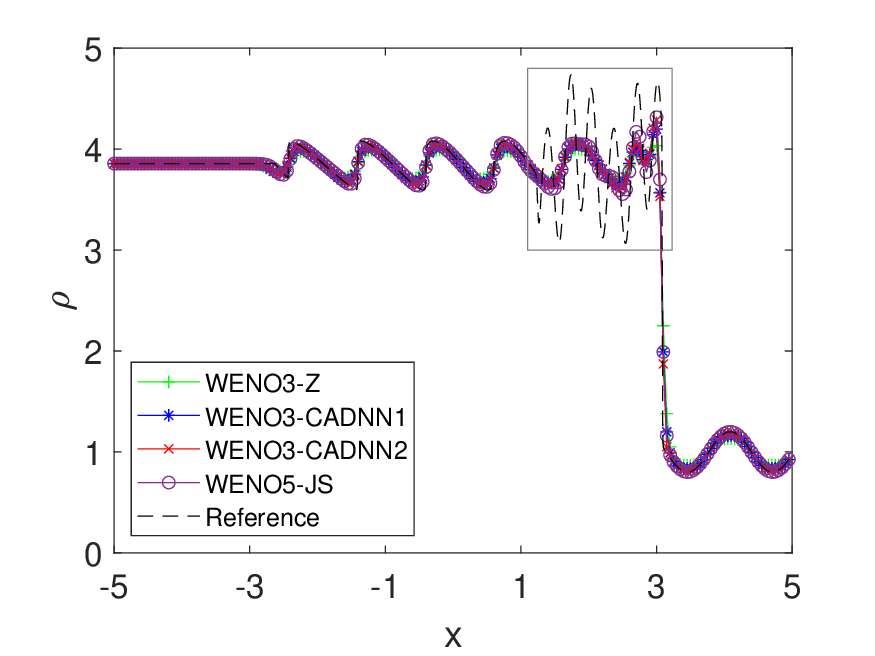}
\includegraphics[height=0.32\textwidth]{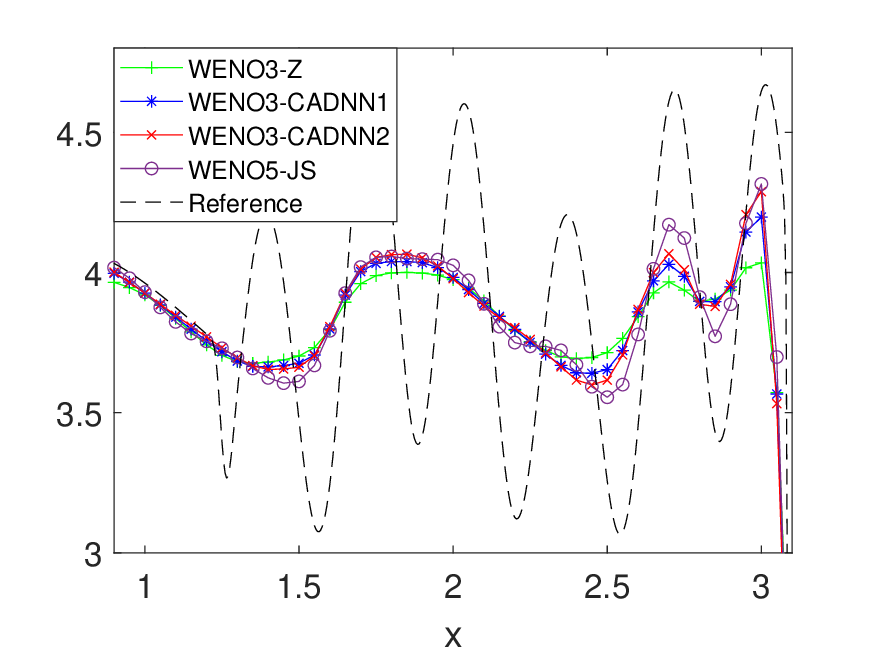}
\caption{Density profiles for Example \ref{ex:shock_entropy_wave} with $k=5$ at $T=2$ (left) and close-up view of the solutions in the box (right) approximated by WENO3-Z (green), WENO3-CADNN1 (blue), WENO3-CADNN2 (red) and WENO5-JS (purple) with $N = 200$. 
The dashed black line is generated by fifth-order WENO5-M with $N = 2000$.}
\label{fig:shock_entropy_wave_k5}
\vspace{-0.15cm}
\end{figure}

\begin{figure}[htbp]
\centering
\includegraphics[height=0.32\textwidth]{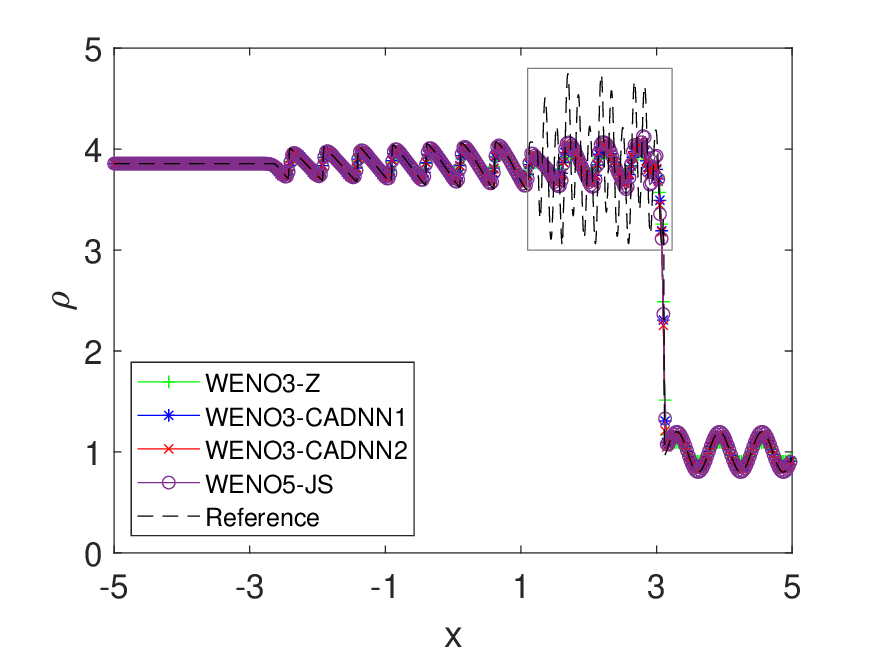}
\includegraphics[height=0.32\textwidth]{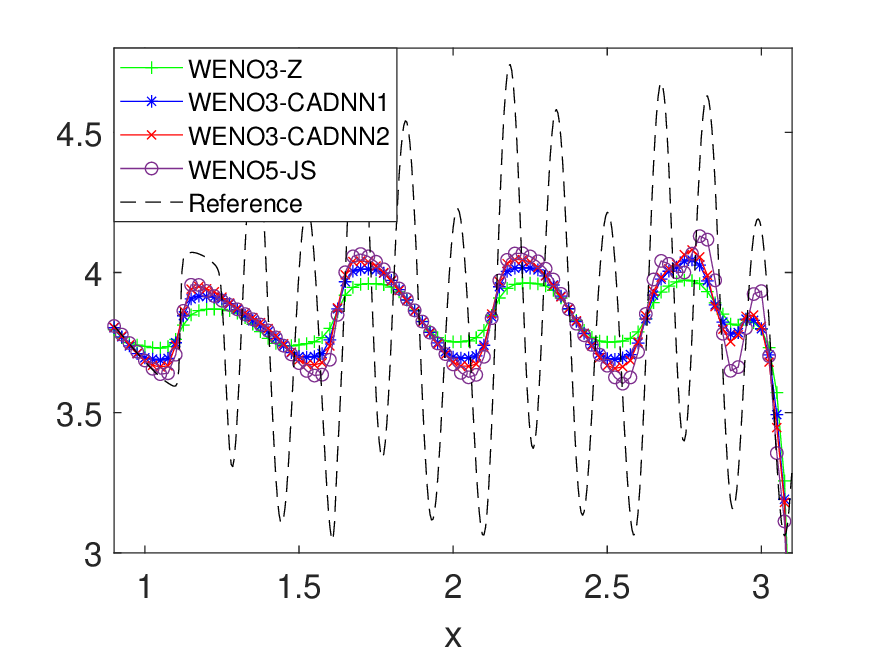}
\vspace{-0.3cm}
\caption{Density profiles for Example \ref{ex:shock_entropy_wave} with $k=10$ at $T=2$ (left) and close-up view of the solutions in the box (right) computed by WENO3-Z (green), WENO3-CADNN1 (blue), WENO3-CADNN2 (red) and WENO5-JS (purple) with $N = 400$.
The dashed black line is generated by fifth-order WENO5-M with $N = 2000$.}
\label{fig:shock_entropy_wave_k10}
\end{figure}

\begin{example} \label{ex:shock_entropy_wave}
The shock entropy wave interaction problem \cite{ShuOsherI}, which involves a Mach 3 shock propagating to the right and interacting with an entropy wave in the density field, is a test case to resolve high-frequency structures. 
The initial condition is defined as
$$
   (\rho, u, P ) = \left\{ 
                    \begin{array}{ll} 
                     (3.857143, \, 2.629369, \, 10.333333), & x < -4, \\ 
                     (1 + 0.2 \sin(kx), \, 0, \, 1),        & x \geqslant -4,
                    \end{array} 
                   \right. 
$$
where $k$ denotes the wave number of the entropy wave.

For $k=5$, the computational domain $[-5, 5]$ is discretized using $N=200$ uniform cells. 
The reference solution is obtained with the fifth-order WENO5-M \cite{Henrick} on a higher resolution of $N = 2000$ points.
Fig. \ref{fig:shock_entropy_wave_k5} plots the numerical solutions of density from WENO3-Z, WENO3-CADNNs and WENO5-JS at $T=2$.

For $k=10$, the same computational domain $[-5, 5]$ is divided uniformly into $N=400$ grid points, and the reference solution is again computed using WENO5-M with $N = 2000$ points.
Fig. \ref{fig:shock_entropy_wave_k10} shows the approximate density profiles by WENO3-Z, WENO3-CADNNs and WENO5-JS at $T=2$.

The numerical results demonstrate that WENO3-CADNNs improve the resolution of the fine-scale entropy waves in comparison with WENO3-Z. 
In particular, due to its lower dissipation, WENO3-CADNN2 achieves accuracy comparable to that of WENO5-JS in some regions, as shown in Figs. \ref{fig:shock_entropy_wave_k5} and \ref{fig:shock_entropy_wave_k10}.
\end{example}

\begin{figure}[htbp]
\centering
\includegraphics[height=0.35\textwidth]{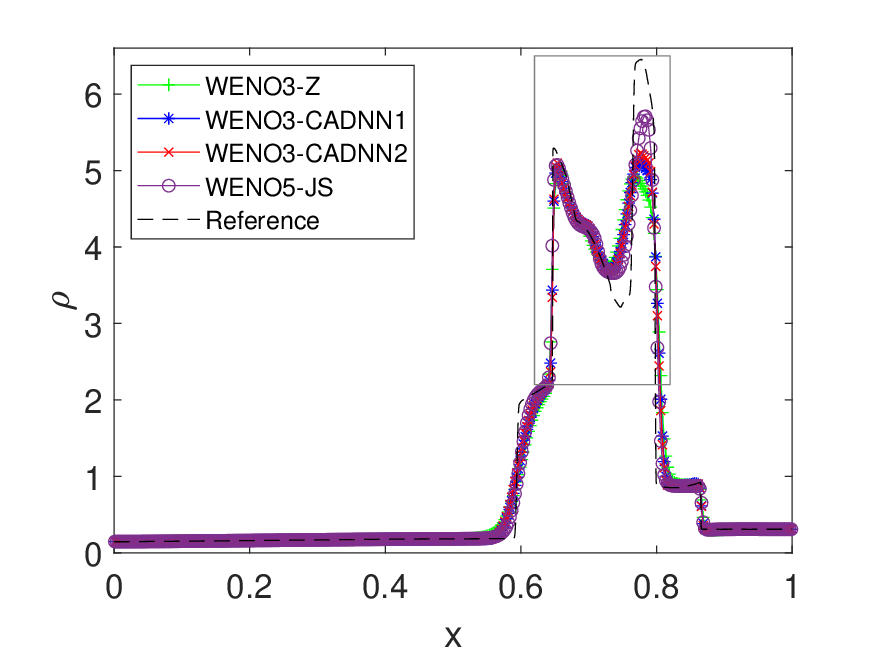}
\includegraphics[height=0.35\textwidth]{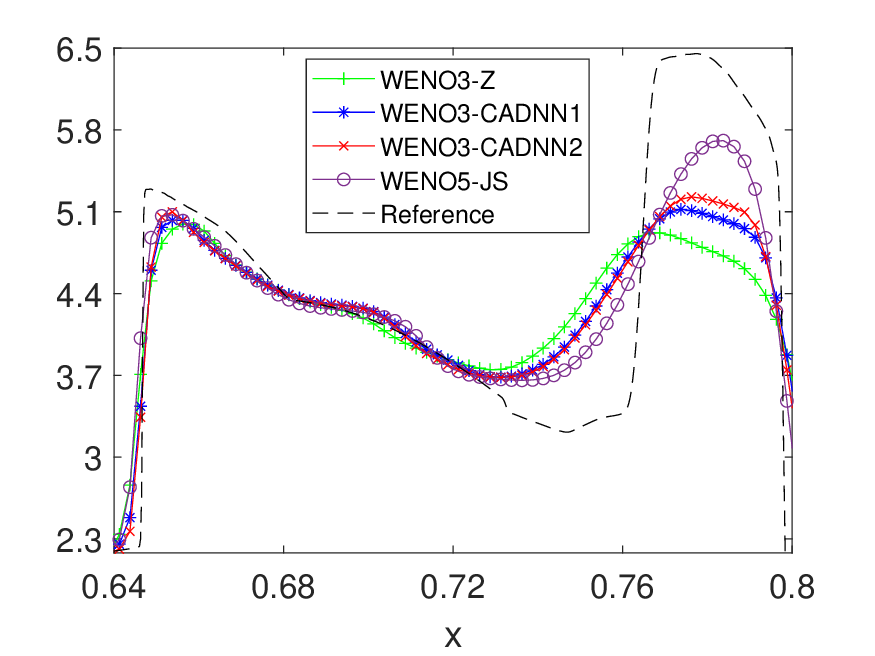}
\caption{Density profiles for Example \ref{ex:interacting_blastwave} at $T=0.038$ (left) and close-up view of the solutions in the box (right) computed by WENO3-Z (green), WENO3-CADNN1 (blue), WENO3-CADNN2 (red) and WENO5-JS (purple) with $N = 400$.
The dashed black lines are generated by fifth-order WENO5-M with $N = 4000$.}
\label{fig:interacting_blastwave}
\end{figure}

\begin{example} \label{ex:interacting_blastwave}
The blast wave interaction problem \cite{Woodward} models the evolution of two blast waves that form, interact, and generate a new contact discontinuity.
The initial condition is given by
$$
   (\rho, u, P ) = \left\{ 
                    \begin{array}{ll} 
                     (1,~0,~1000), & 0 \leqslant x < 0.1, \\ 
                     (1,~0,~0.01), & 0.1 \leqslant x < 0.9, \\
                     (1,~0,~100),  & 0.9 \leqslant x \leqslant 1,
                    \end{array} 
                   \right. 
$$
and reflective boundary conditions are imposed at both ends of the computational domain $[0, 1]$, which is divided into $N = 400$ uniform cells.
We display the numerical density profiles at the final time $T=0.038$ by WENO3-Z, WENO3-CADNNs and WENO5-JS in Fig. \ref{fig:interacting_blastwave}, where the dashed line represents the reference solution calculated from the fifth-order WENO5-M on a refined grid of $N=4000$ points.
It can be seen that WENO3-CADNNs achieve better resolution than WENO3-Z, while WENO5-JS provides superior overall accuracy to the third-order WENO schemes.
\end{example}

\subsection{Two-dimensional system problems}
The two-dimensional Euler equations governing the dynamics of a compressible gas can be expressed in conservative form as
\begin{equation} \label{eq:euler_2d}
 \bfu_t + \bff(\bfu)_x + \bfg(\bfu)_y = 0, 
\end{equation}
where the column vector $\bfu$ of conserved variables, and the flux vectors in the $x$- and $y$-directions, $\bff$ and $\bfg$, are given by
\begin{align*}
      \bfu  &= \left[ \rho, \, \rho u, \, \rho v, \, E \right]^T, \\
 \bff(\bfu) &= \left[ \rho u, \, \rho u^2+P, \, \rho u v, \, u(E+P) \right]^T, \\ 
 \bfg(\bfu) &= \left[ \rho v, \, \rho u v, \, \rho v^2+P, \, v(E+P) \right]^T.
\end{align*}
Here, $\rho$ denotes the density, while $u$ and $v$ represent the velocity components in the $x$- and $y$-directions, respectively. 
The specific kinetic energy $E$ is given by
$$
    E = \frac{P}{\gamma - 1} + \frac{1}{2} \rho (u^2 + v^2),
$$
where $P$ is the pressure and $\gamma = 1.4$ is for the ideal gas, unless otherwise stated.

\begin{figure}[htbp]
\centering
\includegraphics[height=0.28\textwidth]{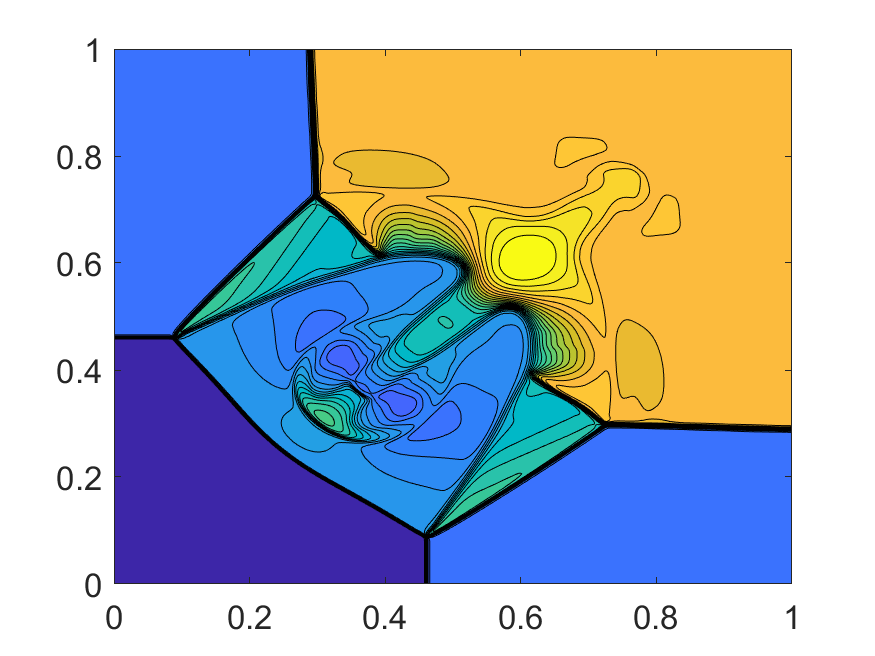} 
\includegraphics[height=0.28\textwidth]{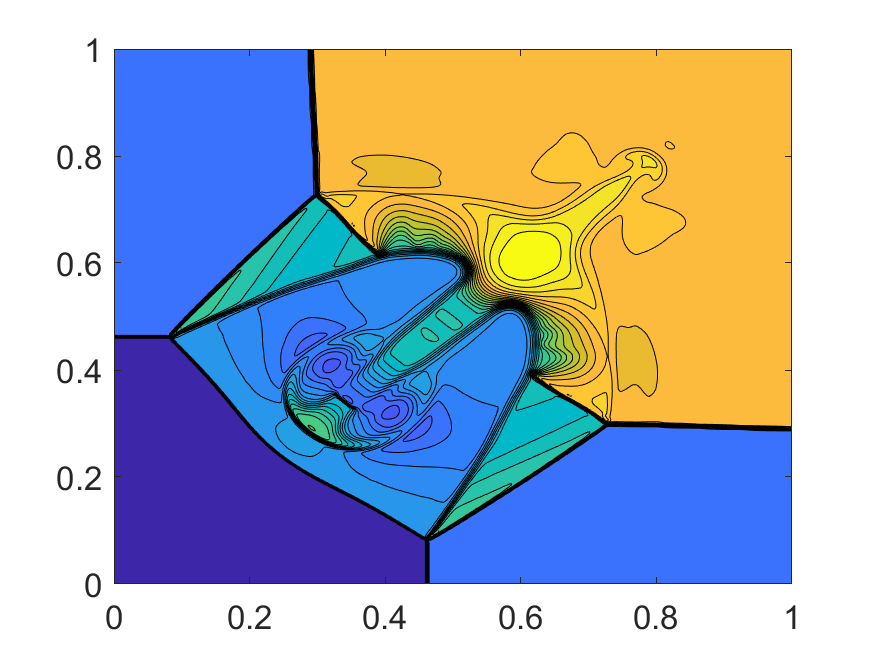}
\includegraphics[height=0.28\textwidth]{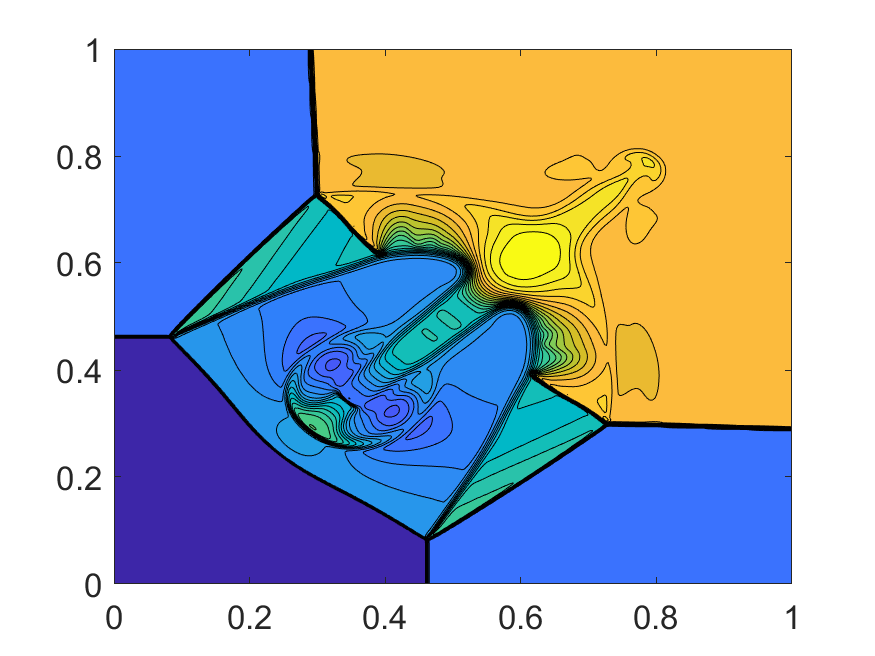}
\includegraphics[height=0.28\textwidth]{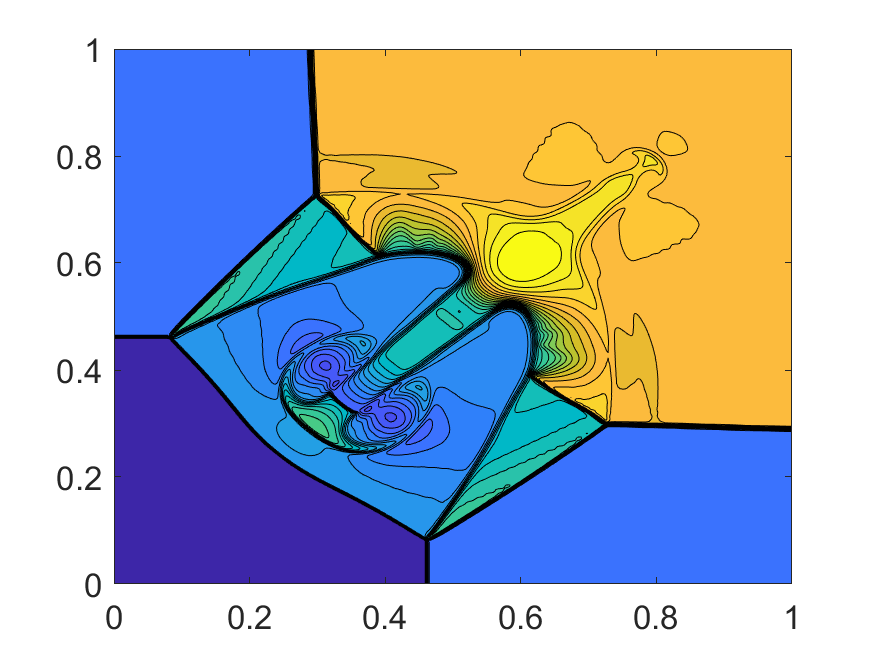}
\vspace{-0.3cm}
\caption{Density in the filled contour plot for Example \ref{ex:euler_2d_riemann} at $T=0.8$ by WENO3-Z (top left), WENO3-CADNN1 (top right), WENO3-CADNN2 (bottom left) and WENO5-JS (bottom right) with $N_x = N_y = 400$.
Each contour plot displays contours at $30$ levels of the density.}
\label{fig:eer2d}
\end{figure}

\begin{example} \label{ex:euler_2d_riemann}
In the first example, we consider the Riemann problem \cite{LaxLiu,Kurganov} for the two-dimensional Euler equations \eqref{eq:euler_2d}, with the following initial condition
$$
   (\rho, u, v, P ) = \left\{ 
                       \begin{array}{ll} 
                        (1.5, \, 0, \, 0, \, 1.5),  & x > 0.8, \, y > 0.8, \\
                        (0.5323, \, 1.206, \, 0, \, 0.3),   & x \leqslant 0.8, \, y > 0.8, \\
                        (0.138, \, 1.206, \, 1.206, \, 0.029),  & x \leqslant 0.8, \, y \leqslant 0.8, \\
                        (0.5323, \, 0, \, 1.206, \, 0.3), & x > 0.8, \, y \leqslant 0.8.
                       \end{array}
                      \right.
$$
We divide the square computational domain $[0, 1] \times [0, 1]$ into $N_x \times N_y = 400 \times 400$ uniform cells.
Fig. \ref{fig:eer2d} presents the numerical solutions of the density at the final time $T=0.8$ by WENO3-Z, WENO3-CADNNs, and WENO5-JS.
It is observed that each scheme captures the reflection shocks and contact discontinuities. 
Additionally, WENO3-CADNN1 and WENO5-JS produce finer roll-up structures of instabilities than WENO3-Z and WENO3-CADNN2.
\end{example}

\begin{figure}[htbp]
\centering
\includegraphics[height=0.3\textwidth]{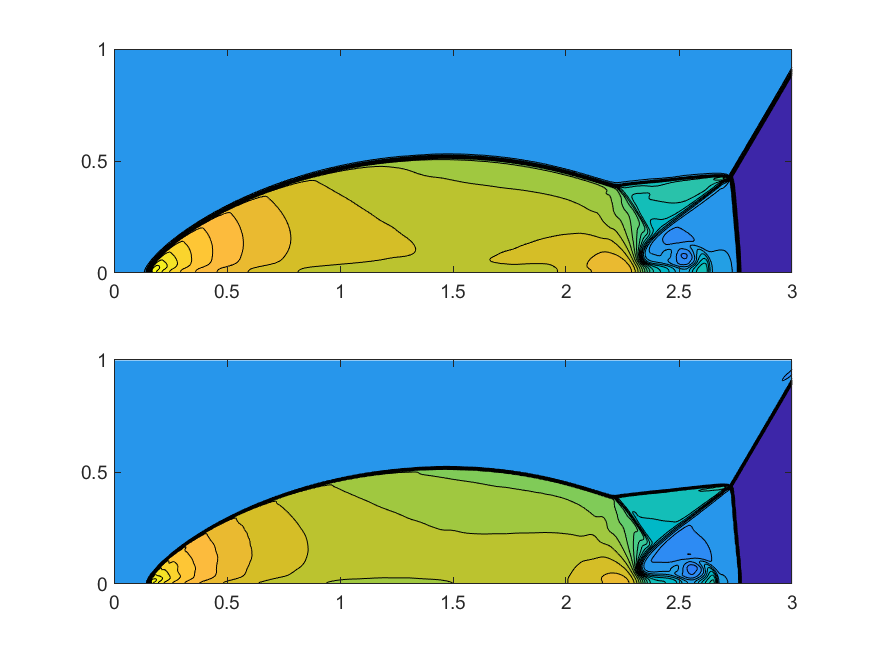}
\includegraphics[height=0.3\textwidth]{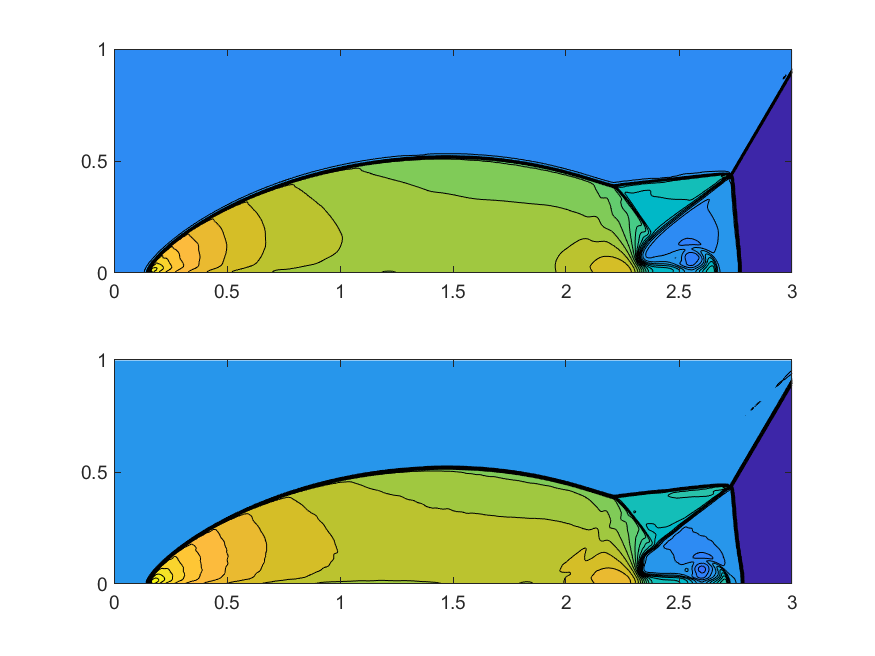}
\vspace{-0.4cm}
\caption{Density in the filled contour plot for Example \ref{ex:double_mach_reflection} at $T=0.2$ by WENO3-Z (top left), WENO3-CADNN1 (top right), WENO3-CADNN2 (bottom left) and WENO5-JS (bottom right) with $N_x = 800$ and $N_y = 200$.
Each contour plot displays contours at $30$ levels of the density.}
\label{fig:dmr}
\end{figure}

\begin{figure}[htbp]
\centering
\includegraphics[height=0.28\textwidth]{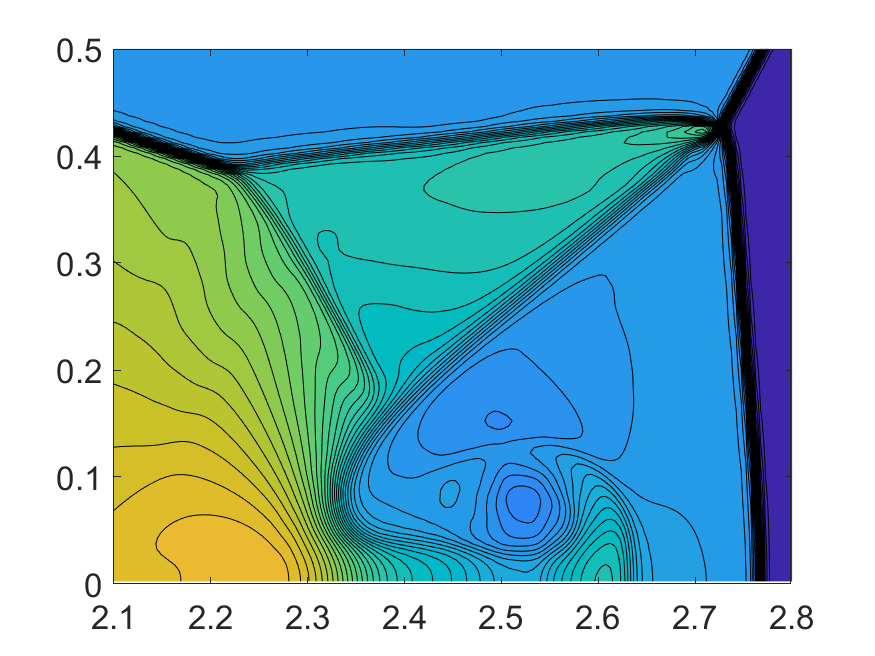}
\includegraphics[height=0.28\textwidth]{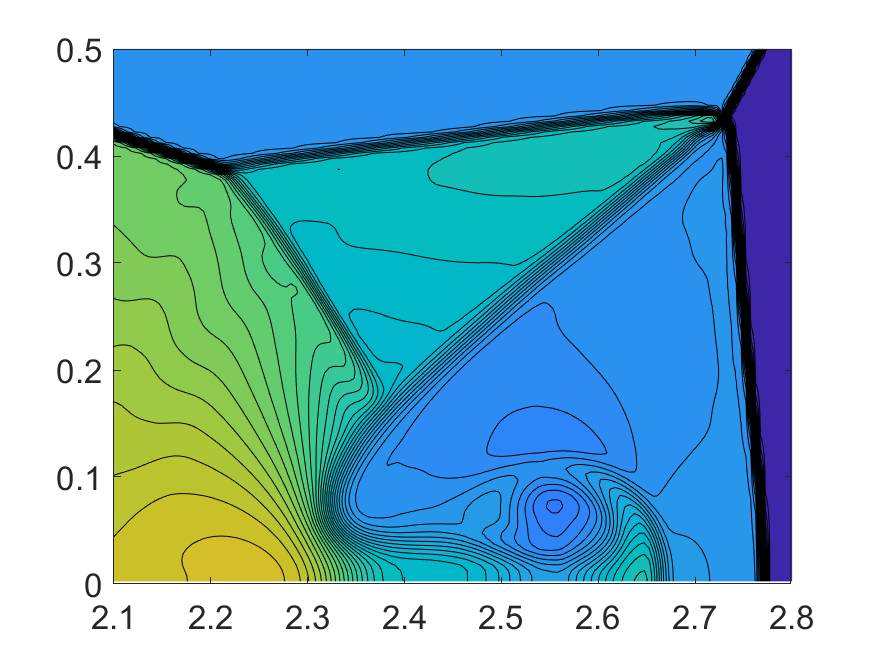}
\includegraphics[height=0.28\textwidth]{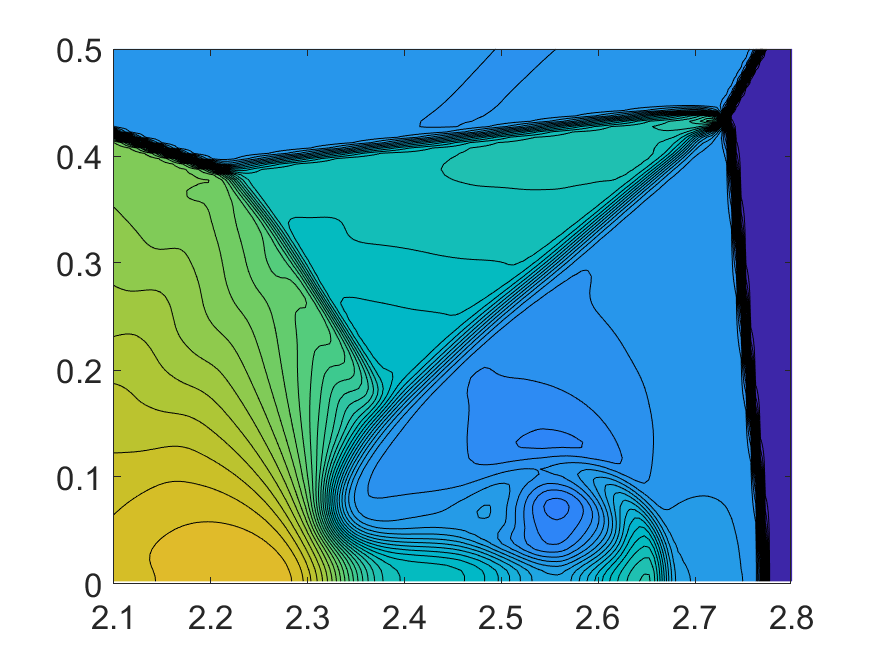}
\includegraphics[height=0.28\textwidth]{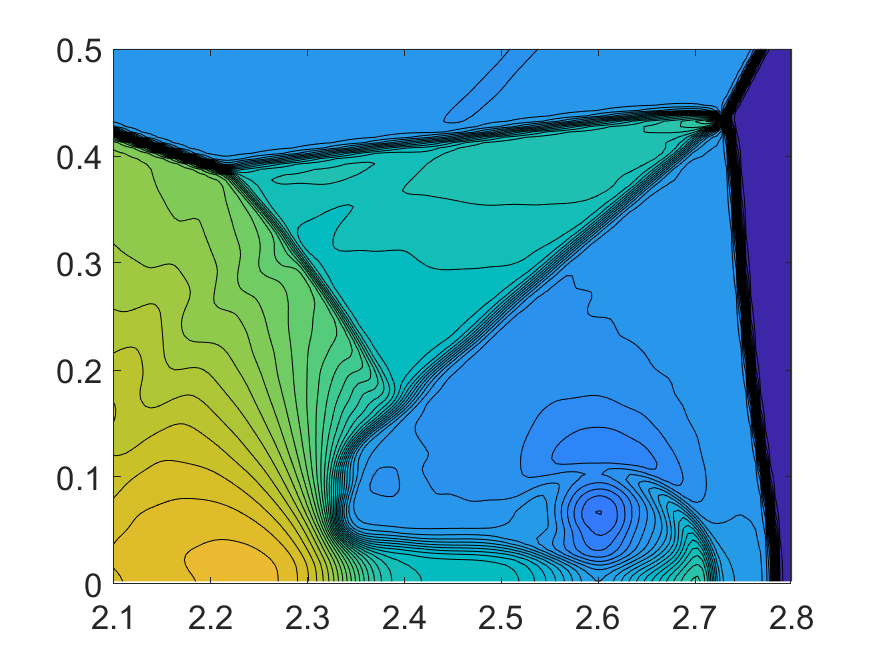}
\caption{Zooming-in density in the filled contour plot for Example \ref{ex:double_mach_reflection} at $T=0.2$ by WENO3-Z (top left), WENO3-CADNN1 (top right), WENO3-CADNN2 (bottom left) and WENO5-JS (bottom right) for the region $[2.2, 2.8] \times [0, 0.5]$.
Each contour plot displays contours at $30$ levels of the density.}
\label{fig:dmr_zoom}
\end{figure}

\begin{example} \label{ex:double_mach_reflection}
We now consider the double Mach reflection problem \cite{Woodward}.
The initial condition is given by 
$$
   (\rho, u, v, P ) = \left\{ 
                       \begin{array}{ll} 
                        (8, \, 8.25 \cos \theta, \, -8.25 \sin \theta, \, 116.5), & x < \frac{1}{6} + \frac{y}{\sqrt{3}}, \\ 
                        (1.4, \, 0, \, 0, \, 1), & x \geqslant \frac{1}{6} + \frac{y}{\sqrt{3}}, 
                       \end{array} 
                      \right. 
$$
with $\theta = \frac{\pi}{6}$.
We divide the computational domain $[0,4] \times [0,1]$ into $N_x \times N_y = 800 \times 200$ uniform cells. 
The simulation is performed up to the final time $T = 0.2$, when the strong shock connecting the contact surface and the transverse wave becomes well-developed.
Fig. \ref{fig:dmr} shows the density approximation by each WENO scheme over the subdomain $[0,3] \times [0,1]$, while Fig. \ref{fig:dmr_zoom} provides a close-up view of the region $[2.2, 2.8] \times [0, 0.5]$.
We can see that WENO3-CADNNs capture the wave interactions more accurately near the second triple point and predict a stronger wall jet than WENO3-Z, with the performance comparable to WENO5-JS.
\end{example}

\begin{figure}[htbp]
\centering
\includegraphics[height=0.35\textwidth]{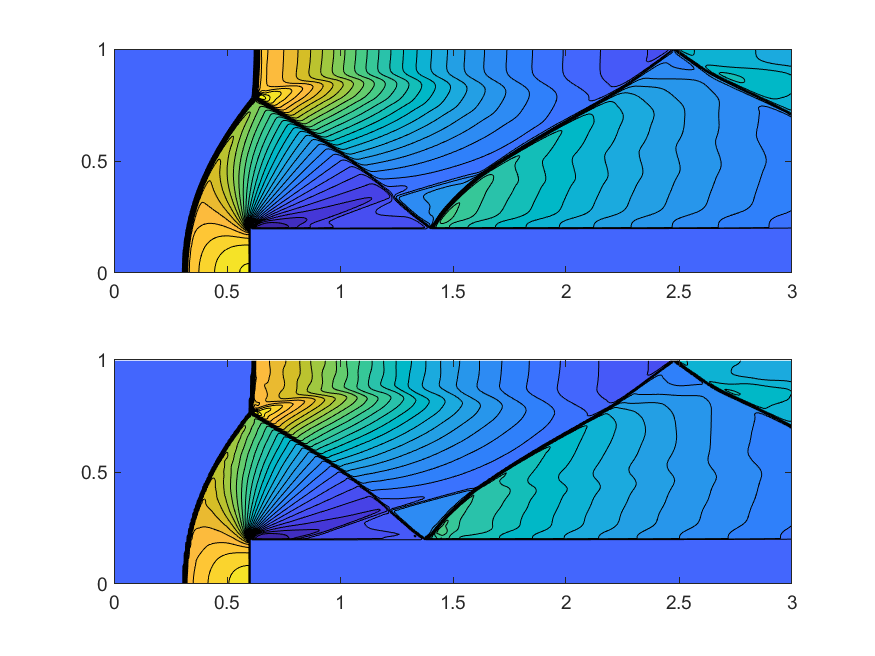}
\includegraphics[height=0.35\textwidth]{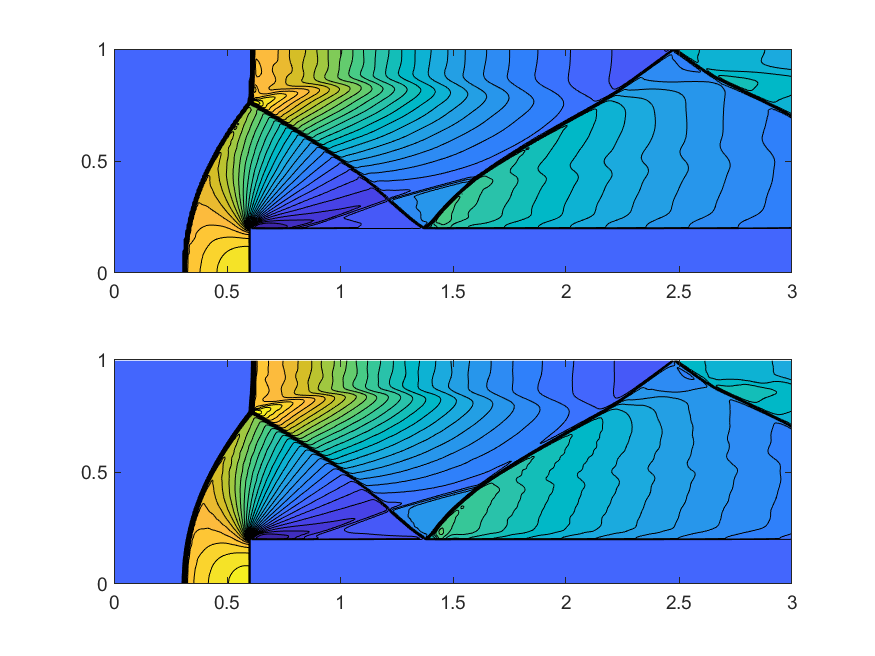}
\caption{Density in the filled contour plot for Example \ref{ex:forward_facing_step} at $T=4$ by WENO3-Z (top left), WENO3-CADNN1 (top right), WENO3-CADNN2 (bottom left) and WENO5-JS (bottom right) with $N_x = 480$ and $N_y = 160$.
Each contour plot displays contours at $30$ levels of the density.}
\label{fig:ffs}
\end{figure}

\begin{example}\label{ex:forward_facing_step}
The forward facing step \cite{Woodward} is a benchmark problem with the initial condition
$$ (\rho, u, v, P ) = (1.4, 3, 0, 1). $$
We divide the computational domain $[0,3] \times [0,1]$ into $N_x \times N_y = 480 \times 160$ uniform grids.
Fig. \ref{fig:ffs} shows the numerical solution on the density at the final time $T=4$.
It is observed that the vortical structures along the slip line by WENO3-CADNNs are less diffused than those by WENO3-Z, and comparable to WENO5-JS.
\end{example}

\begin{figure}[htbp]
\centering
\includegraphics[width=\textwidth]{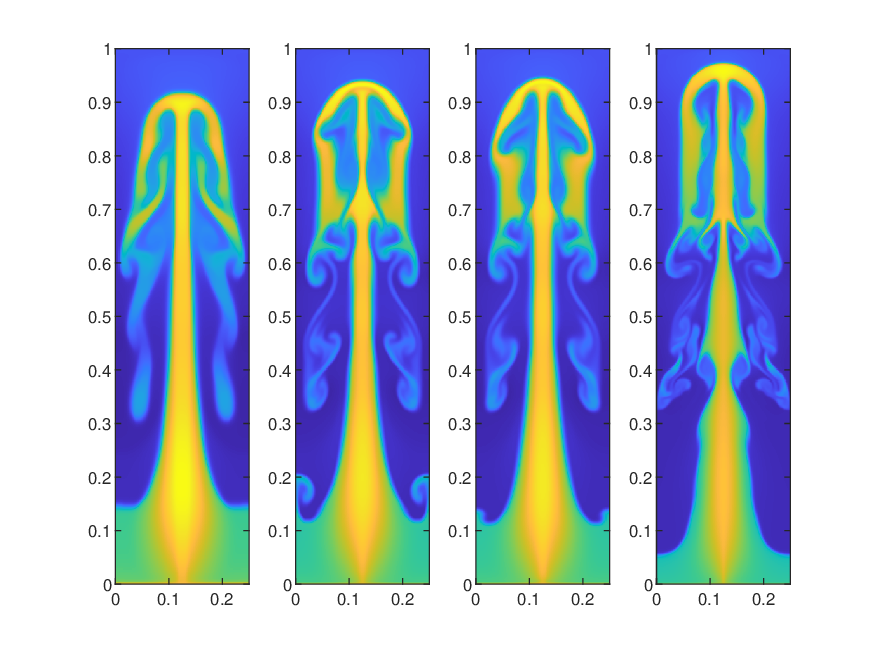}
\caption{Density in the plot of scaled colors for Example \ref{ex:rayleigh_taylor_instability} at $T=2.95$ by WENO3-Z, WENO3-CADNN1, WENO3-CADNN2 and WENO5-JS from left to right with $N_x = 200$ and $N_y = 800$.}
\label{fig:rti}
\vspace{-0.3cm}
\end{figure}

\begin{example}\label{ex:rayleigh_taylor_instability}
The Rayleigh-Taylor instability arises at the interface between two fluids of different densities when a force acts in the direction from the heavy fluid toward the light one.
The instability exhibits a fingering pattern, where bubbles of the light fluid rise into the heavy fluid while spikes of the heavy fluid descend into the light one.
A source term $\bfs(\bfu) = \left[ 0, \, 0, \, \rho, \, \rho v \right]^T$ is added to the right-hand side of the Euler equations \eqref{eq:euler_2d}. 
The initial condition is set as in \cite{ShiZhangShu},
$$
   (\rho, u, v, P ) = \left\{ 
                       \begin{array}{ll} 
                        (2, \, 0, \, -0.025c \cdot \cos(8 \pi x), \, 2y+1), & 0 < y < 0.5, \\ 
                        (1, \, 0, \, -0.025c \cdot \cos(8\pi x), \, y+1.5), & 0.5 \leqslant y \leqslant 1, 
                       \end{array} 
                      \right. 
$$
where $c = \sqrt{\gamma P/\rho}$ is the speed of sound and $\gamma = 5/3$ is the ratio of specific heats.
We compute the numerical solution in the domain $[0, 0.25] \times [0,1]$ with the uniform grids $N_x \times N_y = 200 \times 800$ up to the final time $T=2.95$.
The numerical simulations by WENO3-Z, WENO3-CADNN1, WENO3-CADNN2 and WENO5-JS are plotted in Fig. \ref{fig:rti}, where WENO3-CADNNs resolves finer vortical structures compared to WENO3-Z, while WENO5-JS achieves the best resolution of small-scale vortices because of its higher formal order of accuracy.
\end{example}

\section{Conclusion} \label{sec:conclusions}
In this work, we develop the data-driven WENO weighting function using a feedforward neural network based on the conservative approximation to the derivative.
The network is trained through supervised learning on a labeled dataset created from one-dimensional conservative approximation of the derivative.
The input data consists of four-point stencils sampled from both smooth and discontinuous functions.
The corresponding labels are determined by the regularity of stencils: the exact derivative for the smooth stencil is utilized as the label, whereas we employ the average rate of change for the stencil containing the discontinuity.
To guide the training, we design the loss function comprising three components: the loss of conservative approximation to derivative, the symmetric-balancing term, and the linear term. 
The loss of conservative approximation to derivative minimizes the discrepancy between the approximation to the derivative using the output weights and the label.
The symmetric-balancing term enforces the neural network to satisfy the symmetry property in WENO3-JS and WENO3-Z weighting functions, which suppresses spurious oscillations in the numerical simulations. 
Additionally, we include the linear term with a modified smoothness gauge informed by the smoothness of the stencil, so as to achieve the third-order accuracy in smooth regions and reduce the numerical dissipation near discontinuities.
By tuning the hyperparameters, we obtain two trained models, WENO3-CADNN1 and WENO3-CADNN2, which deliver more accurate solutions than WENO3-Z across all tested cases, and comparable results to WENO5-JS in some instances.

\section*{Acknowledgments}
Jae-Hun Jung is supported by National Research Foundation of Korea (NRF) under the grant number 2021R1A2C3009648, POSTECH Basic Science Research Institute under the NRF grant number 2021R1A6A1A10042944, and partially NRF grant funded by the Korea government (MSIT) (RS-2023-00219980).


\providecommand{\bysame}{\leavevmode\hbox to3em{\hrulefill}\thinspace}
\providecommand{\MR}{\relax\ifhmode\unskip\space\fi MR }
\providecommand{\MRhref}[2]{%
  \href{http://www.ams.org/mathscinet-getitem?mr=#1}{#2}
}
\providecommand{\href}[2]{#2}

\end{document}